\documentclass{article}
\usepackage{amsfonts}
\usepackage{amsmath}
\usepackage{amssymb}
\usepackage{amsxtra}
\usepackage{graphicx}
\usepackage{geometry}
\usepackage{color}
\usepackage{colortbl}
\usepackage{caption}
\usepackage[draft]{varioref}%
\providecommand{\U}[1]{\protect\rule{.1in}{.1in}}
\newtheorem{theorem}{Theorem}

\newtheorem{lemma}[theorem]{Lemma}

\newtheorem{proposition}[theorem]{Proposition}
\newtheorem{remark}[theorem]{Remark}

\numberwithin{equation}{section}
\begin{document}
\title{Asymptotic analysis, polarization matrices and topological derivatives for piezoelectric materials with small voids}

\author{G.Cardone\\University of Sannio - Department of Engineering\\Corso Garibaldi, 107 - 82100 Benevento, Italy\\email: giuseppe.cardone@unisannio.it
\and S.A.Nazarov\\Institute of Mechanical Engineering Problems\\V.O., Bolshoi pr., 61, 199178, St. Petersburg, Russia.\\email: srgnazarov@yahoo.co.uk
\and J.
Sokolowski\\Institut Elie Cartan, UMR 7502 Nancy-Universit\'e-CNRS-INRIA\\ Laboratoire de Math{\'e}matiques, Universit{\'e} Henri Poincar{\'e}\\
Nancy 1, B.P.~239, 54506 Vandoeuvre l{\'e}s Nancy Cedex, France \\email: Jan.Sokolowski@iecn.u-nancy.fr}

\maketitle

\begin{abstract}
Asymptotic formulae for the mechanical and electric fields in a
piezoelectric body with a small void are derived and justified. Such results are
new and useful for applications in the field of design of smart
materials. In this way the topological derivatives of shape
functionals are obtained for piezoelectricity. The asymptotic
formulae are given in terms of the so-called polarization tensors
(matrices) which are determined by the integral characteristics of
voids. The distinguished feature of the piezoelectricity boundary
value problems under considerations is the absence of positive
definiteness of an differential operator which is non
self-adjoint. Two specific Gibbs' functionals of the problem are
defined by the energy and the electric enthalpy. The topological
derivatives are defined in different manners for each of the
governing functionals. Actually, the topological derivative of the
enthalpy functional is local i.e., defined by the pointwise
values of the governing fields, in contrary to the energy
functional and some other suitable shape functionals which admit
non-local topological derivatives, i.e., depending on the whole
problem data. An example with the weak interaction between
mechanical and electric fields provides the explicit asymptotic
expansions and can be directly used in numerical procedures of
optimal design for smart materials.

\medskip

\medskip

\end{abstract}

\section{Introduction}
\label{sec:0}
The paper is devoted to the asymptotic analysis of  boundary value problems for coupled models. The coupling occurs between the mechanical
part which takes the form of the linearized elasticity and governs the stress-strain state of the body, and the electrical part which describes
the electromagnetic field in the body.

From the view point of applications, piezoelectric materials
are of common use in electromechanical sensors and actuators,
e.g., ultrasound transducers in medical imaging and therapy, force
and acceleration sensors, positioning sensors, surface acoustic
wave filters, still with the growing range of applications in
modern technology. Their mode of action is based on the
piezoelectric effect, that couples the electrical and mechanical
behavior of such materials. For the optimal design of
piezoelectric devices, efficient numerical procedures for shape
and topology optimization should be still developed. In the modern
theory of shape optimization it is required that the derivation of
shape and topological derivatives of shape functionals to be
optimized is performed beforehand. From  one side, the derivation
of shape gradients of integral functionals in smooth domains
\cite{sokzol} and non-smooth domains \cite{MaNaPl} (cf.
\cite{na320, na366}) has become a standard procedure. There is no
major difficulty to perform such a shape sensitivity analysis for
the elliptic boundary value problem under considerations. However,
the boundary value problem in piezoelectricity cannot be posed in
such a way that it simultaneously is formally self-adjoint and
possesses a semi-bounded quadratic form. This specific feature
makes the problem more involved from the asymptotic analysis point
of view compared to the pure elasticity or pure electricity
boundary value problems. In addition, the general case of
inhomogeneous and anisotropic body is considered, which also
requires for additional and new technicalities in asymptotic
procedures which is the main subject of the paper. In particular,
different formulations of the piezoelectricity problem (cf.
Sections \ref{2sec:1}, \ref{3sec:1}, \ref{3sec:2}) lead to two
definitions of the polarization matrices which differ one from
another by its properties. Moreover, only the electrical enthalpy,
which is but the governing functional for the piezoelectric media
(see, e.g., \cite{Grin, ikeda, English}) admits the topological
derivative dependent on local chara\-cte\-ristics of mechanical and
electrical fields. Other shape functionals, especially the energy
functional, get the topological derivatives dependent on the
global chara\-cte\-ristics of mechanical and electrical fields. This
acquired trait raises the natural question on the properties of
material derivatives for piezoelectricity in the framework of the
shape sensitivity analysis with smooth or non-smooth boundary
variations, it is clear that the result could be of the same
nature, since the topological derivatives can be identified from
the first order shape gradients by a limit passage e.g. in
elasticity, \cite{sz2001} (cf. also \cite{sz2008}).

In the paper, we restrict ourselves to the asymptotic procedures
of singular domain perturbations which allow us to obtain,  in a
natural way, the topological derivatives of shape functionals. In
principle, the method developed here can be generalized to
characterize the influence on solutions of the non-smooth boundary
variations, therefore, we can derive the shape gradients even in
such a case, e.g., for small defects located close-by the
boundary, including micro-cracks (see \cite{NaSoAA}).

Without entering into details, but with the strong practical
implications in mind, we can claim that some possible applications
of shape optimization in the field concern the design of
electro-acoustic transducers which are constructed with
piezoelectric actuator-patches and capacitative micro-machined
ultrasound transducers. The task for optimal design for a class of
electrostatic-mechanical-acoustic transducers can be e.g., the
topology of electro-acoustic material and the topology of the
electrode-layers, in order to achieve a maximal acoustic pressure,
or a maximal acoustic energy in a specific sub-domains of the
hold-all-domain. We refer the reader e.g. to \cite{English, Grin,
ikeda} for modeling of piezoelectric materials, to
\cite{kaltenbacher} for material tensor identification for such
materials, and to \cite{perla} for control issues.

Our aim is a possible application in shape optimization, thus we
introduce the so-called topological derivatives of shape
functionals for piezoelectric materials. It seems that the models
are not up to now used in applied mathematics for the purposes of
shape optimization, although the smart materials are of common use
in the engineering practice. In shape optimization, the modern
approach to numerical solution, requires the preliminary knowledge
of explicit formulae for shape gradients \cite{sokzol} as well as
of the topological derivatives \cite{sz1997, na320, gagu}. These
formulae are used in the {\it level-set}-type methods which model
the geometrical domain evolution by a zero-level set of solutions
to non-linear hyperbolic equations of the Hamilton-Jacoby type.
The shape gradient are present as the coefficients of the
equations, and the topological derivatives are used to improve the
values of the shape functional under consideration by the
appropriate topology changes, e.g., for the minimization of the
shape functional, the minima of the topological derivative of the
functional indicate the location of a new hole in the geometrical
domain \cite{allaire, fulman1, fulman2}.

\section{Preliminaries.
Problem formulation and
description of results}
\label{sec:1}

\subsection{Shape optimization in piezoelectricity}
\label{1sec:1}

This paper is motivated by the fact that, among numerous
publications on shape optimization, shape sensitivity analysis for
piezoelectric bodies does not exist, although piezoelectric
materials are of extremely wide usage in the modern technologies,
one can think of a simple lighter, available in any supermarket,
or an elaborated computer work-station in a university. One, and
definitely not the only one, distinguishing feature of such smart
materials implies an easy energy transfer in both directions from
mechanical fields to electric fields. The mathematical modeling of
such a phenomenon  leads to serious complications of analysis for
governing partial differential equations  because the corresponding boundary value problem
is not formally self-adjoint in contrast to the boundary value
problems for purely elastic bodies or purely electromagnetic
media. This fact requires for the development of new mathematical
tools and a careful choice of the cardinal shape functional
while neglecting of non-self-adjointness provokes mistakes in
both, mathematical formulae and physical interpretation of the
obtained results (see Remark \ref{r2.99} below).

Introduced in \cite{sz1997}\footnote{Actually, asymptotic formulae
of type \eqref{e0.1} together with the whole asymptotic series for
energy functionals under various singular boundary perturbations
had been derived much earlier in \cite{na83}, although the notion
of the topological derivative is due to \cite{sz1997}.}, the
topological derivative $T(u^{0};\omega_1)$ of a shape functional
$\mathcal{J}$ is intended to describe the change of the functional
$\mathcal{J}$ due to nucleation of holes or voids and allows to
extend possible variations of the shape in an optimization process
\cite{allaire, fulman1, fulman2} in comparison with classical
tools (cf. \cite{sokzol, choi, delfour}),
\begin{equation}\label{e0.1}
\mathcal{J}(u^{h};\Omega(h))=\mathcal{J}(u;\Omega)+h^{\kappa}\mathcal{T}(u;\omega_1)+o(h^{\kappa}),\ h\rightarrow+0,
\end{equation}
In \eqref{e0.1}, $h>0$ is a small parameter, i.e., the diameter of
the opening $\omega_h$ in the entire body
$\Omega\subset\mathbb{R}^{n}$, $u^h$ and $u$ are solutions of the
boundary value problem in
$\Omega(h)=\Omega\setminus\overline{\omega}_h$ and $\Omega$,
respectively, and the exponent $\kappa>0$ depends on the space
dimension $n$ and boundary conditions imposed on the interior
$(n-1)$-dimensional surface $\partial\omega_h$.

Asymptotic analysis of elliptic problems in singularly perturbed
domains, e.g., methods of matched and compound asymptotic
expansions (cf. \cite{Ilin} and \cite{MaNaPl}, respectively), has
become the most appropriate and relevant to obtain {\it almost
explicit} formulae for the topological derivatives as it has been
demonstrated in \cite{na320, na366} and others. We also mention
books \cite{MovMov, ammari} where the subject is studied, to  some
extend, from physical and numerical point of view.

Strangely enough, only self-adjoint problems were heretofore
examined carefully, although  the full-blown approach in
\cite{MaNaPl} can work for general boundary value problems for
elliptic systems. In this paper we partly fill this gap by
adapting formula \eqref{e0.1}   to certain shape functionals for a
piezoelectric body.

The piezoelectricity problems admits two different formulations
with non-symmetric and symmetric but non-semibounded quadratic
forms, the energy and electric enthalpy functionals, respectively.
By means of the Lax-Milgram lemma, the first formulation furnishes
the existence and uniqueness result. At the same time, the
topological derivative of the energy functional is a non-local
characteristics of the piezoelectricity solutions  in
contrast to the pure elasticity problem (see Remark \ref{r2.99}
below), while the absence of this intrinsic property is not caused
by an {\it incorrect} definition \eqref{e0.1}   but occurs as well
for the energy release rate in mechanics of cracks for
piezoelectric media (see Remark \ref{r2.99} again). The fair
explanation, we refer the reader to \cite{Wi} for the modeling
issues, of the latter refers to the electric enthalpy as one of
Gibbs' functional obtained from the energy functional by the
partial Lagrange transform on the  electric components. This is
the electric enthalpy $\mathcal{E}(u^h;\Omega(h))$ (see the definition in
\eqref{e1.16}),    that governs the mechanical electric state of the
piezoelectric body $\Omega(h)$ and, therefore, the second
formulation becomes variational and provides the clear
interpretation of the topo\-lo\-gi\-cal derivative
$\mathcal{T}_{\mathcal{E}}(u;\omega_1)$ in
\begin{equation}\label{e0.2}
\mathcal{E}(u_h;\Omega(h))=\mathcal{E}(u;\Omega)+h^3\mathcal{T}_{\mathcal{E}}(u;\omega_1)+O(h^{4}),\ h\rightarrow+0.
\end{equation}

The indicated peculiarity of the piezoelectricity problem
crucially influences topo\-lo\-gi\-cal derivatives of other shape
functionals, too. For example, the traditional adjoint state (cf.
\cite{choi, sokzol, sz1997}) has to be found out in the formally
adjoint boundary value problem that occasionally underlines its
name.

All the above observations lifts the piezoelectricity problem on
the top of the list of unsolved problems in shape optimization, it
seems that even the classical formulae for material derivatives,
which are not under consideration in the paper, ought to be
revisited.

\subsection{Methods of asymptotic analysis}
\label{2sec:100}

Nowadays there exist several methods to construct asymptotic
expansions of solutions to elliptic boundary value problems in
domains with singular perturbations of boundaries. First of all,
we mention two methods, namely, the method of mathed asymptotic
expansions and the method of compound asymptotic expansions (cf.
monographs \cite{Ilin} and \cite{MaNaPl}, respectively), which in
general appear to be of the same power. Indeed, based on different
asymptotic procedures, the matching procedure and the procedure of
of rearrangement of discrepancies, they result in asymptotic
expansions which differ at the first sight one from another, but
can be readily transformed one into another (we refer to the
introductory chapter 2 in \cite{MaNaPl}). By the way, we silently
use this transformation while presenting at the end of Section
\ref{5sec:2} an alternative way of presentation the asymptotic
form of the derived solution. The method of compound asymptotics
is employed throughout the paper for two reasons. First, the
results given in  \cite{MaNaPl} are obtained in relatively general
formulation which includes the elliptic systems of partial
differential equations not necessarily formally self-adjoint (cf.
discussion in the preceding Section \ref{1sec:1}). On the other
hand, in the monograph \cite{Ilin} the results are established
exclusively for the scalar second-order elliptic equations in the
divergence form. Second, the method of compound asymptotic
expansions is carefully matched with theory of elliptic problems
in domains with conical outlets to infinity, specifically the
exterior domains (cf. [\cite{NaPl}; Ch.6], and \cite{na262}) while
for our purposes this theory is used further to introduce and
investigate the polarization matrices in piezoelectricity.

Since the problem under studies is geometrically very specific, i.e., it concerns only one small opening inside of a domain in $\mathbb{R}^3$, the other methods of asymptotic analysis can be employed. In Remark \ref{r2.11} (2) we mention the case of piecewise constant coefficients which makes suitable an asymptotic analysis of the equivalent boundary integral equations obtained from fundamental solution (cf. \cite{ammari} with similar results in elasticity), although no fundamental matrix is known in piezoelectricity. The other possibilities include among others the homogeneization technique relying on the so-called delute limit (see \cite{JikKozOl}, \cite{Milton}, \cite{Lipton} and many others). However, in our opinion, the method of compound asymptotic expansions is still the most appropriate tool in piezoelectricity in order to investigate asymptotic properties of shape functionals.

\subsection{Constitutive relations in piezoelectricity}
\label{2sec:1}

Let $\Omega\subset\mathbb{R}^{3}$ be a piezoelectric body with the
Lipschitz boundary $\partial\Omega$ and the compact closure
$\overline{\Omega}=\Omega\cup\partial\Omega$. Using the
matrix/column notation (cf. \cite{Lekh, Nabook}), we regard the
displacement vector $u^{\sf{M}}$ as the column
$u^{\sf{M}}=(u_1^{\sf{M}},u_2^{\sf{M}},u_3^{\sf{M}})^\top$ where
$u_{j}^{\sf{M}}$ is the projection of $u$ on the $x_j$-axis of the
fixed Cartesian coordinates system $x=(x_1,x_2,x_3)^{\top}$ and
$\top$ stands for transposition. Together with the electric
potential $u^{\sf{E}}$, the displacements compose the column
$u=(u_1^{\sf{M}},u_2^{\sf{M}},u_3^{\sf{M}},u^{\sf{E}})^{\top}$ of
height 4. The strain column
\begin{equation}\label{e1.1}
\varepsilon^{\sf{M}}(u^{\sf{M}})=(\varepsilon_{11}^{\sf{M}},\varepsilon_{22}^{\sf{M}},\varepsilon_{33}^{\sf{M}},\sqrt{2}\varepsilon_{23}
^{\sf{M}},\sqrt{2}\varepsilon_{31}^{\sf{M}},\sqrt{2}\varepsilon_{12}^{\sf{M}})^\top
\end{equation}
consists of the Cartesian components
$\varepsilon_{jk}^{\sf{M}}=\frac{1}{2}(\partial_{j}u_{k}^{\sf{M}}+\partial_{k}u_{j}^{\sf{M}})$
of the strain tensor and takes the form
$\varepsilon^{\sf{M}}(u^{\sf{M}})=D^{\sf{M}}(\nabla_{x})u^{\sf{M}}$
where
\begin{equation}\label{e1.2}
D^{\sf{M}}(\nabla_x)^\top=\left(
\begin{array}{cccccc}
\partial_1 & 0 & 0 & 0 & 2^{-1/2}\partial_3 & 2^{-1/2}\partial_2 \\
0 & \partial_2 & 0 & 2^{-1/2}\partial_3 & 0 & 2^{-1/2}\partial_1 \\
0 & 0 & \partial_3 & 2^{-1/2}\partial_2 & 2^{-1/2}\partial_1 & 0
\end{array}
\right),
\nabla_x=\left(
\begin{array}{c}
\partial_1 \\
\partial_2 \\
\partial_3
\end{array}
\right),
\partial_j=\frac{\partial}{\partial x_j}.
\end{equation}
We introduce the column
$\varepsilon(u)=(\varepsilon^{\sf{M}}(u^{\sf{M}})^\top,\varepsilon^{\sf{E}}(u^{\sf{E}})^\top)^\top$
where $\varepsilon^{\sf{E}}(u^{\sf{E}})=\nabla_x u^{\sf{E}}$ is
the electric strain column, taken with the sign minus, and
$D(\nabla_x)$ implies a $(9\times4)$-matrix of the first-order
differential operators,
\begin{equation}\label{e1.3}
\varepsilon(u)=D(\nabla_x)u,\ D(\nabla_x)^\top=\left(
\begin{array}{cc}
D^{\sf{M}}(\nabla_x)^\top & \mathbf{0} \\
\mathbf{0}\quad\mathbf{0} & \nabla_x^\top
\end{array}
\right),\
\mathbf{0}=(0,0,0).
\end{equation}

We also assemble  the column $\sigma(u)$ of height $9$ from the
stress column $\sigma^{\sf{E}}(u^{\sf{M}})$ of structure
\eqref{e1.1}   and the electric induction column
$\sigma^{\sf{E}}(u^{\sf{E}})=(\sigma_{1}^{\sf{E}},\sigma_{2}^{\sf{E}},\sigma_{3}^{\sf{E}})^\top$.
In this manner, the constitutive relations of piezoelectricity
(see \cite{Grin, ikeda, English})
\begin{equation}\label{e1.4}
\sigma^{\sf{M}}=A^{\sf{MM}}\varepsilon^{\sf{M}}-A^{\sf{ME}}\varepsilon^{\sf{E}},\ \sigma^{\sf{E}}=A^{\sf{EM}}\varepsilon^{\sf{M}}+A^{\sf{EE}}
\varepsilon^{\sf{E}}
\end{equation}
can be rewritten as follows:
\begin{equation}\label{e1.5}
\sigma(u)=A\varepsilon(u),
\end{equation}
where the matrix $A$ of size $9\times9$,
\begin{equation}\label{e1.6}
A=\left(
\begin{array}{cc}
A^{\sf{MM}} & -A^{\sf{ME}} \\
A^{\sf{EM}} & A^{\sf{EE}}
\end{array}
\right)
\end{equation}
is formed by the symmetric and positive definite $(6\times6)$- and
$(3\times3)$-matrices $A^{\sf MM}$ and $A^{\sf EE}$, respectively
the elastic stiffness matrix and the dielectric permeability
matrix, and the blocks $A^{\sf{ME}}=(A^{\sf{EM}})^\top$ of
piezoelectric moduli. We emphasize that, by its physical nature,
the matrix \eqref{e1.6} is not symmetric provided the
$(6\times3)$-block $A^{\sf{ME}}$ is not null, i.e., the mechanical
and electric fields interact.

The state of the piezoelectric body $\Omega$ is described by the mixed boundary value problem
\begin{align}
D(-\nabla_x)^{\top}A(x)D(\nabla_x)u(x)=f(x),\ x\in\Omega,\label{e1.7} \\
D(n(x))^{\top}A(x)D(\nabla_x)u(x)=g(x),\ x\in\Gamma_{\sigma},\label{e1.8} \\
u(x)=0,\ x\in\Gamma_{u}=\partial\Omega\setminus\overline{\Gamma}_{\sigma},\label{e1.9}
\end{align}
where $n=(n_1,n_2,n_3)^\top$ is the unit vector (column) of the
outward normal. On the right hand-side of the equations \eqref{e1.7}
and \eqref{e1.8}, we have the volume forces
$f^{\sf{M}}=(f_{1}^{\sf{M}},f_{2}^{\sf{M}},f_{3}^{\sf{M}})^\top$
and the surface mechanical loading
$g^{\sf{M}}=(g_{1}^{\sf{M}},g_{2}^{\sf{M}},g_{3}^{\sf{M}})^\top$
together with the volume $f^{\sf{E}}$ and surface $g^{\sf{E}}$
electric charges. The Dirichlet conditions \eqref{e1.9}   mean
that the body is mechanically clamped over the surface $\Gamma_u$
and in contact with an electric conductor. As usually,
$f^{\sf{E}}=0$ and, if the surface $\Gamma_{\sigma}$ is in contact
with a dielectric medium, i.e., vacuum, we are to put
$g_{4}^{E}=0$.

\subsection{Solvability of boundary value problem}
\label{3sec:1}

Let us assume that $mes_2\Gamma_u>0$ and $f\in L^{2}(\Omega)^{4}$, $g\in L^{2}(\Gamma_{\sigma})^{4}$ where $L^{2}(\Xi)$ denote the Lebesque space with the intrinsic inner product $(\ ,\ )_{\Xi}$ and the superscript $4$ indicates the number of components in the vector functions $f$ and $g$.
Notice that the subscript is always omitted in our notation for inner products and norms.

The integral identity (cf. \cite{Lad}) serving for problem \eqref{e1.7}-\eqref{e1.9}, reads as follows:
\begin{equation}\label{e1.10}
Q(u,v;\Omega):=(AD(\nabla_x)u,D(\nabla_x)v)_{\Omega}=(f,v)_{\Omega}+(g,v)_{\Gamma_u},\ v\in\mathring{H}^{1}(\Omega;\Gamma_u)^{4},
\end{equation}
where $\mathring{H}^{1}(\Omega;\Gamma_u)$ denotes the Sobolev space of functions vanishing at $\Gamma_u$.
The left-hand side of \eqref{e1.10}  is understood properly
provided entries of the matrix $A$ are measurable and uniformly
bounded functions in $\Omega$. In addition, for almost all
$x\in\Omega$, we assume the symmetry and positivity properties
\begin{equation}\label{e1.00}
\begin{array}{l}
A^{\sf MM}(x)=A^{\sf MM}(x)^{\top},\ A^{\sf ME}(x)=A^{\sf ME}(x)^{\top},\ A^{\sf EE}(x)=A^{\sf EE}(x)^{\top},\\
c_{\sf M}|a^{\sf M}|^2\leq (a^{\sf M})^\top A^{\sf MM}(x)\leq C_{\sf M}|a^{\sf M}|^2,\ a^{\sf M}\in\mathbb{R}^6,\\
c_{\sf E}|a^{\sf E}|^2\leq (a^{\sf E})^\top A^{\sf EE}(x)\leq
C_{\sf E}|a^{\sf E}|^2,\ a^{\sf E}\in\mathbb{R}^3,
\end{array}
\end{equation}
where $c_{\sf M}$, $C_{\sf M}$ and $c_{\sf E}$, $C_{\sf E}$ are
positive constants. We emphasize that no positivity restriction is
imposed on the piezoelectric moduli in $A^{\sf ME}$.

Although in the case $A^{\sf{ME}}\neq0$ the sesquilinear form
$Q(\cdot,\cdot;\Omega)$ cannot be an inner product on the Hilbert
space $\mathring{H}^{1}(\Omega;\Gamma_u)^{4}$ due to the {\it
wrong} sign on $A^{\sf{ME}}$ in \eqref{e1.6}, the Lax-Milgram
lemma ensures the following assertion because of the formula
\begin{equation}\label{e1.000}
Q(u,u;\Omega)=(A^{\sf MM}D^{\sf M}(\nabla_x)u^{\sf M},D^{\sf M}(\nabla_x)u^{\sf M})_\Omega+(A^{\sf EE}\nabla_x u^{\sf E},\nabla_x u^{\sf E})_\Omega
\geq c\|u;H^1(\Omega)\|^2
\end{equation}
caused by the Poincar\'e inequality for $u^{\sf E}$ and the Korn
inequality for $u^{\sf M}$ (see \cite{DuLi, KoOl} and others).

\begin{proposition}\label{p1.1}
Under  the conditions  \eqref{e1.00},  \eqref{e1.000}, the problem
\eqref{e1.10} admits a unique solution
$u\in\mathring{H}^{1}(\Omega;\Gamma_u)^{4}$, and the following
estimate is valid:
\begin{equation}\label{e1.11}
\|u;H^{1}(\Omega)\|\leq c_\Omega(\|f;L^{2}(\Omega)\|+\|g;L^{2}(\Gamma_\sigma)\|).
\end{equation}
\end{proposition}
Unfortunately, the problem \eqref{e1.10}   is non variational.
Indeed, the energy functional $\mathcal{U}$,
\begin{align}
\mathcal{U}(u;\Omega)=\frac{1}{2}(AD(\nabla_x)u,D(\nabla_x)u)_\Omega-\mathcal{A}(u;\Omega), \label{e1.12} \\
\mathcal{A}(u;\Omega):=(f,u)_\Omega+(g,u)_{\Gamma_\sigma},\label{e1.13}
\end{align}
is but the sum of the mechanical and electric energy functionals
\begin{align}
\mathcal{U}^{\sf{M}}(u^{\sf{M}};\Omega)=\frac{1}{2}(A^{\sf{MM}}D^{\sf{M}}(\nabla_x)u^{\sf{M}},D^{\sf M}(\nabla_x)u^{\sf M})
_\Omega-(f^{\sf{M}},u^{\sf{M}})_\Omega-(g^{\sf{M}},u^{\sf{M}})_{\Gamma_\sigma},\label{e1.14}\\
\mathcal{U}^{\sf{E}}(u^{\sf E};\Omega)=\frac{1}{2}(A^{\sf{EE}}\nabla_x u^{\sf{E}},\nabla_x u^{\sf{E}})
_\Omega-(f^{\sf{E}}, u^{\sf{E}})_\Omega-(g^{\sf E}, u^{\sf{E}})_{\Gamma_{\sigma}},\label{e1.15}
\end{align}
while a stationary point of \eqref{e1.12}   becomes a solution of
the problem \eqref{e1.10}   with the block-diagonal
$(9\times9)$-matrix $diag\{A^{\sf{MM}},A^{\sf{EE}}\}$, i.e., the
variational problem does not accept an interaction of the mechanical
and electric fields (cf. an example in Section \ref{sec:examples}).

It is known (see, e.g., \cite{Wi}) that the electric enthalpy
$\mathcal{E}$,
\begin{align}
\mathcal{E}(u;\Omega)=\frac{1}{2}(A_{(-)}D(\nabla_x)u,D(\nabla_x)u)_{\Omega}-\mathcal{R}(u;\Omega),\label{e1.16}\\
\mathcal{R}(u;\Omega)=(f^{\sf{M}},u^{\sf{M}})_{\Omega}+(g^{\sf{M}},u^{\sf{M}})_{\Gamma_\sigma}-(f_{4}^{\sf{E}},u_{4}^{\sf{E}})
_\Omega-(g_{4}^{\sf{E}},u_{4}^{\sf{E}})_{\Gamma_\sigma},\label{e1.17}
\end{align}
gives rise to the variational formulation of the piezoelectricity problem
\begin{equation}\label{e1.18}
Q_{(-)}(u,v;\Omega):=(A_{(-)}D(\nabla_x)u,D(\nabla_x)v)_\Omega=\mathcal{R}(v;\Omega),\ v\in\mathring{H}^{1}(\Omega;\Gamma_\sigma)^{4},
\end{equation}
where the matrix $A_{(-)}$ is composed from blocks in \eqref{e1.6}   as follows
\begin{equation}\label{e1.19}
A_{(-)}=\left(
\begin{array}{cc}
A^{\sf{MM}} & A^{\sf{ME}} \\
A^{\sf{EM}} & -A^{\sf{EE}}
\end{array}
\right)
\end{equation}
The matrix \eqref{e1.19}, in contrast to the matrix $A$, is
symmetric, however, neither  matrix  \eqref{e1.19},   nor  the
quadratic form on the left-hand side of \eqref{e1.18}  is
positive definite. Thus, a solution
$u\in\mathring{H}^{1}(\Omega;\Gamma_\sigma)^{4}$ is a stationary
point of the functional \eqref{e1.16}   but $u$ cannot be any
minimizer of the electric enthalpy $\mathcal{E}(u;\Omega)$.

The integral identity \eqref{e1.18} with the test function
$v_{(-)}=(v_{1}^{\sf{M}},v_{2}^{\sf{M}},v_{3}^{\sf{M}},-v^{\sf{E}})$
trans\-forms into the problem \eqref{e1.10}. The inverse
transformation is also available. These facts prove that the
problem \eqref{e1.18} inherits the unique solvability from
\eqref{e1.10} in Proposition \ref{p1.1}.
\begin{remark}\label{r1.2}
The integral identity is formally obtained by the multiplying system \eqref{e1.7}   with $v$ scalarly and integrating by parts.
Using $v_{(-)}$ as the multiplier, one arrives at \eqref{e1.18}.
This explains the equivalency of the problems.
\end{remark}
\begin{remark}\label{r1.2bis}
In Section \ref{3sec:2} we shall outline a different way to modify the piezoelectricity problem in order to study properties of the mechanical and electric fields on the base of known results. This approach is related to studies \cite{Milton}, \cite{CherGib}, and others, on variational formulations of elliptic problems describing processes in media with complex material coefficients (phase changes due to material properties). In fact we need this technique only for one reason, to maintain the so-called polynomial property \cite{na217, na226}, thus the results
\cite{Milton}, \cite{CherGib}
are not applied in the paper. We recall that the polynomial property allows to describe all required attributes of the exterior boundary value problem of piezoelectricity by simple algebraic calculations (cf. review \cite{na262}).
\end{remark}

The electric enthalpy is but the difference of  elastic energy \eqref{e1.14}   and  electric energy \eqref{e1.15}.
Expression \eqref{e1.13}   implies the external work.
Being the difference of the mechanical and electric external works, the component $\mathcal{R}(u;\Omega)$ of the electric enthalpy has no physical meaning as a whole.
Nevertheless, in Section \ref{2sec:3} we shall observe that asymptotic formulae for $\mathcal{E}(u;\Omega)$ become meaningful while the analogous formulae for $\mathcal{U}(u;\Omega)$ look rather queer.

\subsection{Structure of the paper}
\label{4sec:1}

In Section \ref{sec:2} the asymptotic analysis of
the piezoelectricity problem for the body $\Omega(h)$ with a small
void $\overline{\omega_h}$ is performed (see \eqref{e2.1}). The
applied here asymptotic procedure \cite[Ch.4]{MaNaPl} requires for
introduction of an intrinsic integral characteristics of the void
$\overline{\omega_1}$ in the homogeneous piezoelectric space
$\mathbb{R}^3$,  the polarization matrix $M(A^0,\omega)$ of
size $9\times 9$ (see formulae \eqref{e2.34}-\eqref{e2.36}).
Theorem \ref{t2.8} establishes general properties of the
polarization matrix, see also \eqref{e5.21} for the case of weak
interaction between mechanical and electric fields. The
polarization matrix appears in the asymptotic expansion of the
boundary layer term at infinity that also permits in Section
\ref{5sec:2} to complete the asymptotic ansatz of the solution to
the piezoelectricity problem in $\Omega(h)$. The asymptotics
constructed in Section \ref{sec:2} is justified in Section
\ref{1sec:3}. In Section \ref{2sec:3} the asymptotics of the
energy and electric enthalpy functionals are analysed, while in
Section \ref{3sec:3} rather arbitrary shape functional is
considered and the corresponding adjoint state is detected. The
paper is completed by inquiring into a piezoelectric body with a
weak interaction of the mechanical and electric fields. All
asymptotic formulae derived in the paper are made more explicit in
such a case due to the fact that for pure electricity and pure
elasticity the polarization matrices are known explicitly for many
canonical shapes (see, respectively, \cite{PoSe}, \cite{na100,
pos2, ammari} and others).
\section{Asymptotic analysis}
\label{sec:2}

\subsection{The problem with an interior singular perturbation in the domain}
\label{1sec:2}

Let $\omega$ be an open set in $\mathbb{R}^{3}$ with a Lipschitz boundary and a compact closure.
We assume that both $\Omega$ and $\omega$ contain the coordinate origin $\mathcal{O}$.
Given a small dimensionless parameter $h\in(0,h_0]$, we introduce the sets
\begin{equation}\label{e2.1}
\omega_h=\{x:\xi:=h^{-1}x\in\omega\},\ \Omega(h)=\Omega\setminus\overline{\omega}_h.
\end{equation}
The bound $h_0>0$ is chosen such that $\overline{\omega}_h\subset\Omega$ for $h\in(0,h_0]$.
By rescaling, we reduce a characteristic size of $\Omega$ and $\omega$ to the unit and make the coordinates $x$ and $\xi$ dimensionless.

Supposing $\Omega(h)$ to be a connected set, we consider the piezoelectricity problem in the domain $\Omega(h)$, namely,
\begin{align}
D(-\nabla_x)^{\top}A(x)D(\nabla_x)u^h(x)=f(x),\ x\in\Omega(h),\label{e2.2}\\
D(n(x))^{\top}A(x)D(\nabla_x)u^h(x)=g(x),\ x\in\Gamma_\sigma,\label{e2.3}\\
D(n^{h}(x))^{\top}A(x)D(\nabla_x)u^h(x)=0,\ x\in\partial\omega_h,\label{e2.4}\\
u^h(x)=0,\ x\in\Gamma_u.\label{e2.5}
\end{align}
In \eqref{e2.4}, $n^h$ stands for the outward normal on
$\partial\omega_h$. Since the Neumann conditions are imposed on
the boundary of $\omega_h$, there is no traction on
$\partial\omega_h$ and the opening $\overline{\omega_h}$ is filled
with a dielectric medium. This problem, of course, ought to be
reformulated as either  integral identity \eqref{e1.10}, or
\eqref{e1.18}   in the function space
$\mathring{H}^{1}(\Omega(h);\Gamma_u)^{4}$, hence
\begin{equation}\label{e2.555}
Q(u^h,v^h;\Omega(h))=(f,v^h)_{\Omega(h)}+(g,v^h)_{\Gamma_\sigma},\ v^h\in\mathring{H}^1(\Omega(h);\Gamma_u)^4.
\end{equation}
Proposition \ref{p1.1} remains valid for the problem \eqref{e2.555}   in the domain $\Omega(h)$.

For $h=0$, the opening $\omega_h$ disappears and the singularly
perturbed problem \eqref{e2.2}-\eqref{e2.5}   becomes the original
problem \eqref{e1.7}-\eqref{e1.9}. In order to describe the
behavior of the solution
$u^h\in\mathring{H}^{1}(\Omega(h);\Gamma_u)^4$ as
$h\rightarrow+0$, we have to assume an additional smoothness of
the matrix $A$, for example, in the ball $\mathbb{B}_R=\{x:|x|<R\}$ the
inclusion
\begin{equation}\label{e2.6}
A\in C^{2,\alpha}(\overline{\mathbb{B}}_R)^{9\times9}
\end{equation}
is valid, where $C^{k,\alpha}(\Xi)$ is the H\"older space with the standard norm
\begin{equation*}
\|v;C^{k,\alpha}(\Xi)\|=\sum_{j=1}^{k}\sup_{x\in\Xi}|\nabla_{x}^{j}v(x)|+\sup_{x,y\in\Xi}|x-y|^{-\alpha}|\nabla_{x}^{k}v(x)-\nabla_{y}^{k}v(y)|
\end{equation*}
and $\nabla_{x}^{k}v$ denotes the family of all order $k$ derivatives of $v$.
Since the matrix differential operator
\begin{equation}\label{e2.7}
L(x,\nabla_x)=D(-\nabla_x)^{\top}A(x)D(\nabla_x)
\end{equation}
is elliptic (see Section \ref{3sec:2} below), a solution $u\in
H^1(\mathbb{B}_R)^4$ of system \eqref{e1.7}   in
$\mathbb{B}_R$ with the right-hand side
\begin{equation}\label{e2.8}
f\in C^{0,\alpha}(\mathbb{B}_R)^4,\ \alpha\in(1/2,1),
\end{equation}
falls into the space $C^{2,\alpha}(\mathbb{B}_{R^{\prime}})^4$ for any $R^{\prime}\in(0,R)$.
This fact is due to local estimates of solutions to elliptic systems \cite{ADN2}.
Note that \eqref{e2.8}   provides the estimate
\begin{equation}\label{e2.9}
|f(x)-f(0)|\leq c|x|^{\alpha},\ x\in\mathbb{B}_R.
\end{equation}

We also need the Taylor formula
\begin{equation}\label{e2.10}
|u(x)-d(x)a-D(x)^{\top}\varepsilon^0-U(x)|\leq c|x|^{2+\alpha},\ x\in\mathbb{B}_{R^{\prime}},
\end{equation}
where $D(x)^{\top}$ is the matrix in \eqref{e1.3}   under the
substitution $\nabla_x\mapsto x$,
\begin{equation}\label{e2.11}
\varepsilon^0=D(\nabla_x)u(0)\in\mathbb{R}^9,
\end{equation}
$d(x)a$ with $a\in\mathbb{R}^7$ implies a rigid motion in the mechanical component and a constant potential in the electric one,
\begin{equation}\label{e2.12}
d(x)=\left(
\begin{array}{cc}
d^{\sf{M}}(x) & 0 \\
\mathbf{0}\quad\mathbf{0} & 1
\end{array}
\right),\ d^{\sf{M}}(x)=\left(
\begin{array}{cccccc}
1 & 0 & 0 & 0 & -2^{-1/2}x_3 & 2^{-1/2}x_2 \\
0 & 1 & 0 & 2^{-1/2}x_3 & 0 & -2^{-1/2}x_1 \\
0 & 0 & 1 & -2^{-1/2}x_2 & 2^{-1/2}x_1 & 0
\end{array}
\right).
\end{equation}
We emphasize a similarity of the matrices $D^{\sf{M}}(x)^{\top}$ and $d^{\sf{M}}(x)$.
Finally, $U$ in \eqref{e2.10}   is a quadratic term, i.e.,
\begin{equation}\label{e2.00}
U(tx)=t^2U(x),\ t>0,\ x\in\mathbb{R}^3.
\end{equation}
\begin{remark}\label{r2.1}
The factor $\sqrt{2}$ is present in  the strain column \eqref{e1.1} in
order to equalize the natural norms for tensors of rank $2$ with
the norms of corresponding columns of height $6$. As a result, an
orthogonal transformation of the Cartesian coordinate system $x$
implies the orthogonal transformations for all columns introduced
to replace tensors (see, e.g., \cite[Ch.2]{Nabook}). By the factor
$2^{-1/2}$ in \eqref{e2.12}, we also achieve the relations
\begin{equation}\label{e2.13}
\begin{array}{l}
D(\nabla_x)D(x)^{\top}=\mathbb{I}_{9\times9},\ D(\nabla_x)d(x)=\mathbb{O}_{9\times7},\\
d(\nabla_x)^{\top}d(x)|_{x=0}=\mathbb{I}_{7\times7},\ d(\nabla_x)^{\top}D(x)^{\top}|_{d=0}=\mathbb{O}_{7\times9},
\end{array}
\end{equation}
where $\mathbb{I}_{n\times n}$ and $\mathbb{O}_{m\times n}$ stand for the unit and null matrices of size $n\times n$ and $m\times n$, respectively.
Notice that \eqref{e2.11}   follows from the first couple of the relations \eqref{e2.13}   and our way to write the Taylor formula.
\end{remark}

By \eqref{e2.6}, we particularly obtain
\begin{equation}\label{e2.14}
A(x)=A^0+\sum_{j=1}^{3}x_jA^j+\widetilde{A}(x),\ |\widetilde{A}_{pq}(x)|\leq c|x|^2,\ x\in\mathbb{B}_R,
\end{equation}
with the constant $(9\times9)$-matrices $A^j$ so that matrix \eqref{e2.7}   of differential operator gets the decomposition
\begin{equation}\label{e2.15}
L(x,\nabla_x)=L^0(\nabla_x)+L^{\prime}(x,\nabla_x)+\widetilde{L}(x,\nabla_x).
\end{equation}
Inserting the Taylor formula for $u$ into the equation \eqref{e1.7}   and using \eqref{e2.14}   yield
\begin{equation}\label{e2.16}
L^0(\nabla_x)U(x)-\sum_{j=1}^{3}D(e_j)^{\top}A^j\varepsilon^0=f(0).
\end{equation}
Here $e_j=(\delta_{j,1},\delta_{j,2},\delta_{j,3})^{\top}$.
Since $U$ is quadratic in $x$ (see \eqref{e2.00}), the first term on the left hand-side is independent of $x$.
\begin{remark}\label{r2.10}
To guarantee formulae \eqref{e2.9}   and \eqref{e2.10}   with $\alpha\in(0,1/2)$, we could assume $f\in H^2(\mathbb{B}_R)^4$ while deriving
 $u\in H^4(\mathbb{B}_{R^{\prime}})^4$ from local estimates for solutions of elliptic systems (see \cite{ADN2}).
This is due to the Sobolev embedding theorem $H^{l+2}\subset C^{l,\alpha}$ in $\mathbb{R}^3$ for
any $\alpha\in(0,1/2)$.
However, in Theorem \ref{t3.1} and Remark \ref{r3.2} we shall see that we really need $\alpha>1/2$.
The latter requires, for example, $f\in H^3(\mathbb{B}_R)^4$, and, therefore, we prefer here to use the H\"older scale.
\end{remark}

\subsection{The asymptotic ansatz}
\label{2sec:2}

Based on general results in \cite{MaNaPl} on the asymptotic structure of solutions to elliptic boundary value problems in a domain with singular perturbations of the boundary, we accept the following
asymptotic ansatz for the solution $u^h$ of
problem \eqref{e2.2}-\eqref{e2.5}  :
\begin{equation}\label{e2.17}
u^h(x)=u(x)+\chi(x)(hw^1(\xi)+h^2w^2(\xi))+h^3\mathbf{u}(x)+\dots
\end{equation}
Here $u$ is a solution of the limit problem \eqref{e1.7}-\eqref{e1.9},
$w^1$ and $w^2$ are terms of the boundary layer type, and
$\mathbf{u}$ is the main regular corrector.
The boundary layer terms are treated in Sections \ref{3sec:2} and \ref{4sec:2}, respectively, and the regular corrector in Section \ref{5sec:2} below.
The cut-off function
$\chi\in C_{c}^{\infty}(\Omega)$ is equal to one in the ball
$\mathbb{B}_{R/3}$ and null outside $\mathbb{B}_{2R/3}$ so that, now,  we
fix $h_0>0$ such that $\omega_h\subset\mathbb{B}_{R/3}$ for
$h\in(0,h_0]$.

\begin{remark}
\label{RemR1}
The boundary layer solutions $w^1$ and $w^2$ are constructed in Sections \ref{3sec:2} and \ref{4sec:2}, respectively, along with their decompositions at infinity. In Section \ref{5sec:2} the main terms of the decompositions, compose the right-hand side of a problem of type \eqref{e1.7}-\eqref{e1.9} for the regular solution $\mathbf{u}$.
\end{remark}

In view of \eqref{e2.1}, the coordinate dilation
$x\mapsto\xi=h^{-1}x$ removes the boundary $\partial\Omega$ close
to infinity and the formal limit passage $h\rightarrow+0$ makes
the exterior domain $\Xi=\mathbb{R}^3\setminus\overline{\omega}$
from the nucleated domain $\Omega(h)$. Moreover, the decomposition
\eqref{e2.15}   yields
\begin{equation}\label{e2.18}
L(x,\nabla_x)=L(h\xi,h^{-1}\nabla_\xi)=h^{-2}L^0(\nabla_\xi)+h^{-1}L^{\prime}(\xi,\nabla_\xi)+\dots
\end{equation}
Similarly, for the Neumann boundary operator $N^h(x,\nabla_x)$ on the left hand-side of \eqref{e2.4}, we have
\begin{equation}\label{e2.19}
N^h(x,\nabla_x)=h^{-1}N^0(\xi,\nabla_\xi)+h^0N^{\prime}(\xi,\nabla_\xi)+\dots
\end{equation}
where
\begin{equation}\label{e2.20}
N^0(\xi,\nabla_\xi)=D(n^{\omega}(\xi))^{\top}A^0D(\nabla_x),\ N^{\prime}(\xi,\nabla_\xi)=D(n^\omega(\xi))^{\top}\sum_{j=1}^{3}\xi_jA^jD(\nabla_\xi),
\end{equation}
and $n^\omega$ is the unit vector of the outward normal on $\partial\omega$.

Let us derive the exterior boundary value problems for $w^1$ and
$w^2$. First, we insert the ansatz \eqref{e2.17}   into
\eqref{e2.2}, make use of the expansion \eqref{e2.19}, and collect
coefficients written in the fast variables $\xi$ for similar
powers of the small parameter $h$. As a result, we obtain systems
of differential equations in $\Xi$ for $w^1$ and $w^2$ (see
\eqref{e2.21}   and \eqref{e2.22}   below). Second, we calculate
the discrepancy left by the leading asymptotic term $u(x)$ in the
boundary conditions \eqref{e2.4}. Namely, by means of
\eqref{e2.10}, \eqref{e2.19}, we derive that
\begin{equation*}
N^h(x,\nabla_x)u(x)=D(n^{\omega}(\xi))^{\top}\left(A^0+h\sum_{j=1}^{3}\xi_jA^j\right)\varepsilon^0+hN^0(\xi,\nabla_\xi)U(\xi)+\dots
\end{equation*}
Finally, we write  the problems
\begin{equation}\label{e2.21}
\begin{array}{l}
L^0(\nabla_\xi)w^1(\xi)=0,\ \xi\in\Xi,\\
N^0(\xi,\nabla_\xi)w^1(\xi)=-D(n^\omega(\xi))^{\top}A^0\varepsilon^0,\ \xi\in\partial\omega,
\end{array}
\end{equation}
and
\begin{equation}\label{e2.22}
\begin{array}{l}
L^0(\nabla_\xi)w^2(\xi)=-L^{\prime}(\xi,\nabla_\xi)w^1(\xi),\ \xi\in\Xi,\\
N^0(\xi,\nabla_\xi)w^2(\xi)=-N^\prime(\xi,\nabla_\xi)w^1(\xi)-N^\prime(\xi,\nabla_\xi)D(\xi)^{\top}\varepsilon_0-N^0(\xi,\nabla_\xi)U(\xi),\ \xi\in\partial\omega.
\end{array}
\end{equation}

\subsection{The exterior problem in piezoelectricity}
\label{3sec:2}

The polynomial property  \cite{na217, na262} of a formally
self-adjoint system of differential equations delivers plenty of
results for the exterior boundary value problem in $\Xi$ such as
the ellipticity, the solvability, asymptotic expansions of
solutions, and intrinsic integral characteristics, i.e. the
polarization matrices (see \cite[Ch.6]{NaPl}, \cite{na262, na239}
and \cite{na320} in shape optimization). As it has been mentioned,
the piezoelectricity system \eqref{e1.7}   is not formally
self-adjoint, however,  introducing the imaginary potential
$iu_{4}^{\sf{E}}$ (see \cite[Example 1.13]{na262}) and the column
$u_{(i)}=(u_{1}^{\sf{M}},u_{2}^{\sf{M}},u_{3}^{\sf{M}},iu_{4}^{\sf{E}})^{\top}$
brings the sesquilinear form
\begin{equation}\label{e2.23}
q_{(i)}(u_{(i)},v_{(i)};\Xi)=(A_{(i)}^{0}D(\nabla_\xi)u_{(i)},D(\nabla_\xi)v_{(i)})_{\Xi}
\end{equation}
where $i$ is the imaginary unit and $A_{(i)}^{0}$ stands for  modified matrix \eqref{e1.19},
\begin{equation}\label{e2.24}
A_{(i)}^{0}=\left(
\begin{array}{cc}
A^{0\sf{MM}}, & iA^{0\sf{ME}}\\
iA^{0\sf{EM}}, & A^{0\sf{EE}}
\end{array}
\right)=A_{(Re)}^{0}+iA_{(Im)}^{0},
\end{equation}
while both $A_{(Re)}^{0}$ and $A_{(Im)}^{0}$ are real symmetric
and $A_{(Re)}^{0}$ is positive definite. The sesquilinear form
\eqref{e2.23} is not Hermitian in the case
$A^{0\sf{ME}}\neq\mathbb{O}_{6\times3}$, but it enjoys the
polynomial property \cite{na217, na226, na262}:
\begin{equation}\label{e2.25}
q_{(i)}(u_{(i)},u_{(i)};\Upsilon)=0\ \Longleftrightarrow\ u_{(i)}\in\mathcal{P}|_{\Upsilon},
\end{equation}
where $\Upsilon$ is any domain in $\mathbb{R}^3$ and
$\mathcal{P}=\{p:p(x)=d(x)a,\ a\in\mathbb{C}^7\}$ is a polynomial
subspace of dimension $7$ generated by the matrix in
\eqref{e2.12}.

The above observations made in \cite{na226, na262} and the
investigation scheme \cite[Ch.6]{NaPl} provide all results we
formulate below with exception for the polarization matrix and
here  the most attention is paid to this integral characteristics
of the opening $\overline{\omega}$ in the homogeneous
piezoelectric space.

Let $V_{0}^{1}(\Xi)$ be the Kondratiev space \cite{Ko} obtained by
the completion of the linear space
$C_{c}^{\infty}(\overline{\Xi})$ (infinitely differentiable
functions with compact supports) with respect to the Dirichlet
integral norm $\|\nabla_\xi w;L^2(\Xi)\|$. Applying the
one-dimensional Hardy inequality in the radial variable
$\rho=|\xi|$, we use the equivalent norm
\begin{equation}\label{e2.265}
\|w;V_{0}^{1}(\Xi)\|=(\|\nabla_\xi w;L^{2}(\Xi)\|^2+\|\rho^{-1}w;L^2(\Xi)\|^2)^{1/2}.
\end{equation}

The problem \eqref{e2.21}  with the right-hand side $g\in
L^2(\partial\omega)^4$ in the Neumann boundary conditions can be
reformulated as the integral identity, similarly to \eqref{e1.10}
\begin{equation}\label{e2.26}
(A^0D(\nabla_\xi)w,D(\nabla_\xi)v)_{\Xi}=(g,v)_{\partial\omega},\ v\in V_{0}^{1}(\Xi)^4.
\end{equation}

\begin{proposition}\label{p2.3}
For any $g\in L^{2}(\partial\omega)^4$, the problem \eqref{e2.26}
has a unique solution $w\in V_{0}^{1}(\Xi)^4$ and the estimate
$\|w;V_{0}^{1}(\Xi)\|\leq c\|g;L^2(\partial\omega)\|$ is valid.
\end{proposition}

Although $\partial\omega$ and $g$ are not smooth, the solution $w$
in Proposition \ref{p2.3} is infinitely differentiable outside of
any neighborhood $\mathcal{V}$ of the set $\overline{\omega}$
(recall the local estimates in \cite{ADN2} mentioned above). To
describe the behavior of $w(\xi)$ as $\rho\rightarrow\infty$, we
introduce the fundamental matrix $\Phi(x)$ of size $4\times4$ for
the operator $L^0(\nabla_\xi)$ in $\mathbb{R}^3$ (see
\cite{GeShi,Herm}). This matrix is positive homogeneous of degree
$-1$, namely,
\begin{equation}\label{e2.27}
\Phi(t\xi)=t^{-1}\Phi(\xi),\ t>0,\ \xi\in\mathbb{R}^3\setminus\{0\}.
\end{equation}

The next assertion is due to \cite{Ko}, \cite{MaPl2} (see also \cite{Pazy} and,
e.g., \cite[Ch.6]{NaPl}).

\begin{proposition}\label{p2.4}
The solution $w\in V_{0}^{1}(\Xi)^4$ of the problem \eqref{e2.26}   admits the asymptotic form
\begin{align}
w(\xi)=(d(-\nabla_\xi)^{\top}\Phi(\xi)^{\top})^{\top}a+(D(-\nabla_\xi)\Phi(\xi)^{\top})^{\top}b+\widetilde{w}(\xi),\label{e2.28}\\
|\nabla_{\xi}^{k}\widetilde{w}(\xi)|\leq c_k\rho^{-3-k},\ k\in\mathbb{N}_0=\{0,1,2,\dots\},\ \xi\in\mathbb{R}^3\setminus\mathcal{V},\label{e2.288}
\end{align}
where $a\in\mathbb{R}^7$ and $b\in\mathbb{R}^9$ while $|a|+|b|\leq
c\|g;L^2(\partial\omega)\|$.
\end{proposition}
\begin{remark}\label{RemR2}
Formulae \eqref{e2.28}-\eqref{e2.288} can be derived from the integral representation of the solution $w$ through the fundamental matrix \eqref{e2.27}. In this way decomposition \eqref{e2.288} is obtained from the Taylor formula in inverted variables $\xi |\xi|^{-2}$. Observe that the columns of matrices $D(x)^\top$ and $d(x)$ in \eqref{e1.3} and \eqref{e2.12} form a basis in the linear subspace of dimension $16$ of columns linearly dependent on variables $x=(x_1, x_2, x_3)$. We emphasize that the matrix notation of elasticity relations combined with the polynomial property allow us to write the complete decomposition of $w(\xi)$ for $|\xi|\to\infty$ in a condensed and convenient form for further applications.
\end{remark}
\begin{remark}\label{r2.5}
Formula \eqref{e2.28}   contains the matrices $d$ and $D$ in \eqref{e2.12}   and \eqref{e1.3}.
Let $d^1(\xi),\dots,d^7(\xi)$ be columns of $d(\xi)$ and let $D_1(\xi),\dots,D_9(\xi)$ be strings of $D(\xi)$.
Then we rewrite \eqref{e2.28}   in the form of strings
\begin{equation*}
w(\xi)^{\top}=\sum_{j=1}^{7}a_jd^j(-\nabla_\xi)^{\top}\Phi(\xi)^{\top}+\sum_{k=1}^{9}b_kD_k(-\nabla_\xi)\Phi(\xi)^{\top}+\widetilde{w}(\xi)^{\top}.
\end{equation*}
Therefore, the asymptotic terms detached in \eqref{e2.28}  are but
a linear combination of columns of the fundamental matrix
$\Phi(\xi)$ (with the coefficients $a_1$, $a_2$, $a_3$ and $a_7$;
cf.\eqref{e2.12}) and of the first-order derivatives of the
columns (with the coefficients $a_4$, $a_5$, $a_6$ and
$b_1,\dots,b_9$).
\end{remark}

The columns $d^1,\dots,d^7$ satisfy  the homogeneous problem
\eqref{e2.21}. However, the columns are not in  the weighted space
$V_{0}^{1}(\Xi)^4$ by the lack of their  decay rate and, hence,
$d^j(\xi)$ are not solutions of the homogeneous $(g=0)$ problem
\eqref{e2.26} in Proposition \ref{p2.4}. According to the general
method \cite{MaPl1} such solutions are used to compute the
coefficients in the asymptotic expansion \eqref{e2.28}. We are
going to use this method twice. First, we observe that the
right-hand side $g$ in \eqref{e2.21}   verifies the orthogonality
conditions
\begin{equation}\label{e2.29}
\int\limits_{\partial\omega}d(\xi)^{\top}g(\xi)ds_\xi=0\in\mathbb{R}^7.
\end{equation}
Indeed, by \eqref{e2.13}, we get
\begin{equation}\label{e2.30}
\left(\int\limits_{\partial\omega}d(\xi)^{\top}g(\xi)ds_\xi\right)^{\top}=-(A^0\varepsilon^0)^{\top}\int\limits_{\partial\omega}D(n^\omega(\xi))
d(\xi)ds_\xi=-(A^0\varepsilon^0)^{\top}\int\limits_\omega D(\nabla_\xi)d(\xi)d\xi=0.
\end{equation}
\begin{proposition}\label{p2.6}
Under  orthogonality condition \eqref{e2.29}, the column $a\in\mathbb{R}^7$ in \eqref{e2.28}   vanishes.
\end{proposition}

The proof is commented in Remark \ref{r2.0}.

Let $W^j\in V_{0}^{1}(\xi)^4$ be a solution to the problem
\eqref{e2.26}   with the specific right-hand side
\begin{equation}\label{e2.31}
g^j(\xi)=-D(n^\omega(\xi))^{\top}A^0\mathbf{e}_j;
\end{equation}
here $j=1,\dots,9$,
$\mathbf{e}_j=(\delta_{j,1},\dots,\delta_{j,9})^{\top}$ is the
unit column in $\mathbb{R}^9$, and $\delta_{j,k}$ stands for the
Kronecker symbol. Recalling the problem \eqref{e2.21}   for the
boundary layer term $w^1$, we see that
\begin{equation}\label{e2.32}
w^1(\xi)=W(\xi)\varepsilon^0
\end{equation}
with the $(4\times9)$-matrix function $W$ composed from the columns $W^1,\dots,W^9$ of height $4$,
\begin{equation}\label{e2.33}
W=(W^1,\dots,W^9).
\end{equation}

By Proposition \ref{p2.6} and the relation \eqref{e2.30}, we conclude the expansions
\begin{equation}\label{e2.34}
W^j(\xi)^{\top}=\sum_{p=1}^{9}M_{jp}D_p(\nabla_\xi)\Phi(\xi)^{\top}+\widetilde{W}^j(\xi)^{\top}
\end{equation}
where the remainders $\widetilde{W}^j(\xi)$ obey the estimates
\eqref{e2.288}. The coefficients $M_{jp}$ in \eqref{e2.34}   form
the matrix of size $9\times9$
\begin{equation}\label{e2.35}
M=M(A^0,\omega)
\end{equation}
which, in the analogy with \cite{na100, na280, NazSokSPN} and
others, is called {\it the polarization matrix} of the opening
$\omega$ in the homogeneous piezoelectric space.

As in Section \ref{3sec:1}, our study of general properties of
\eqref{e2.35}   relies on both formulations \eqref{e1.10}   and
\eqref{e1.18}   of the piezoelectricity problem. Hence, we have to
perform the same sign changes as in \eqref{e1.19},
\begin{equation}\label{e2.36}
M=\left(
\begin{array}{cc}
M^{\sf{MM}} & M^{\sf{ME}} \\
M^{\sf{EM}} & M^{\sf{EE}}
\end{array}
\right)\mapsto M_{(=)}=\left(
\begin{array}{cc}
M^{\sf{MM}} & -M^{\sf{ME}}\\
M^{\sf{EM}} & -M^{\sf{EE}}
\end{array}
\right).
\end{equation}
\begin{theorem}\label{t2.8}
Entries of the modified polarization matrix $M_{(=)}$ satisfy the relation
\begin{equation}\label{e2.355}
(M_{(=)})_{jp}=-Q_{(-)}^{0}(W^j,W^p;\Xi)-(A_{(-)}^{0})_{jp}mes_3\omega,\ j,p=1,\dots,9,
\end{equation}
where $Q_{(-)}^{0}$ is the quadratic form in \eqref{e1.18}   with the matrix $A_{(-)}^{0}=A_{(-)}(0)$ (see \eqref{e1.19}   and \eqref{e2.14}).
\end{theorem}

{\bf Proof.} By \eqref{e2.31}   and \eqref{e2.13}, the sum
$\mathcal{W}^j(\xi)=D_j(\xi)^{\top}+W^j(\xi)$ verifies the
homogeneous problem \eqref{e2.21}. In the method \cite{MaPl1}
these solutions play the same role as it was registered for the
columns $d^1,\dots,d^7$ above Proposition \ref{p2.6}. We
underline that the vector function
\begin{equation}\label{e2.365}
\mathcal{W}_{(-)}^{j}=(\mathcal{W}_{1}^{j\sf{M}},\mathcal{W}_{2}^{j\sf{M}},\mathcal{W}_{3}^{j\sf{M}},-\mathcal{W}^{j\sf{E}})^{\top}
\end{equation}
verifies a homogeneous boundary value problem which is formally
adjoint for \eqref{e2.21}   and involves the differential
operators $L_{(\top)}^{0}$ and $N_{(\top)}^{0}$ constructed  from
$L^0$ and $N^0$ in \eqref{e2.15}   and \eqref{e2.20},
respectively, by replacing $A^0$ with the transposed matrix
$(A^0)^{\top}$. Clearly,
$L_{(\top)}^{0}(\nabla_\xi)=L^0(\nabla_\xi)^{\ast}$ is the
formally adjoint for the differential operator $L^0(\nabla_\xi)$.

We insert $W^j$ and $\mathcal{W}_{(-)}^{p}$ into the Green formula
written for the truncated domain $\Xi_R=\Xi\cap\mathbb{B}_R$ and
choose the radius of the ball $\mathbb{B}_R=\{\xi:|\xi|<R\}$ such
that the sphere $\mathbb{S}_R=\partial\mathbb{B}_R$ envelopes the
set $\overline{\omega}$. We have
\begin{equation}\label{e2.375}
(L^0W^j,\mathcal{W}_{(-)}^{p})_{\Xi_R}+(N^0W^j,\mathcal{W}_{(-)}^{p})_{\partial\omega\cup\mathbb{S}_R}=(W^j,L_{(\top)}^{0}
\mathcal{W}_{(-)}^{p})_{\Xi_R}+(W^j,N_{(\top)}^{0}\mathcal{W}_{(-)}^{p})_{\partial\omega\cup\mathbb{S}_R}.
\end{equation}
Since $L^0\mathcal{W}^j=0$ provides $L_{(\top)}^0\mathcal{W}_{(-)}^{j}=0$, the integrals over $\Xi_R$ in \eqref{e2.375}   vanish.
Furthermore, $N_{(\top)}^{0}(\xi,\nabla_\xi)\mathcal{W}_{(-)}^{p}(\xi)=0$, $\xi\in\partial\omega$.
Thus, \eqref{e2.375}   converts into
\begin{equation}\label{e2.37}
(N^0W^j,\mathcal{W}_{(-)}^{p})_{\partial\omega}=(W^j,N_{(\top)}^{0}\mathcal{W}_{(-)}^{p})_{\mathbb{S}_R}-(N^0W^j,\mathcal{W}_{(-)}^{p})_{\mathbb{S}_R}
\end{equation}
where $N^0(\xi,\nabla_\xi)=D(|\xi|^{-1}\xi)^{\top}A^0D(\xi)$ on the sphere $\mathbb{S}_R$.

Taking into account the estimates \eqref{e2.288}  for
$\widetilde{W}^j$ and the concomitant estimates
$|\nabla_\xi^kW^p(\xi)|\leq c_p\rho^{-1-k}$, we obtain that the
right-hand side $I_{right}^{jp}$ of \eqref{e2.375}   satisfies
\begin{equation*}
I_{right}^{jp}=(\Sigma^j,N_{(\top)}^{0}D_{p(-)}^\top)_{\mathbb{S}_R}+O(R^{-1})
\end{equation*}
where $\Sigma^j$ means the asymptotic term detached  in
\eqref{e2.34}   and $D_{p(-)}(\xi)^{\top}$ is a column of the
matrix $D(\xi)^{\top}$ transformed according to \eqref{e2.31}.
Understanding integrals over the ball $\mathbb{B}_R$ in the
framework of the  theory of distributions and using the Green formula, we obtain
\begin{eqnarray}
\nonumber I_{right}^{jp} & = & (L^0\Sigma^j,D_{p(-)}^\top)_{\mathbb{B}_R}-(\Sigma^j,L_{(\top)}^0D_{p(-)}^\top)_{\mathbb{B}_{R}}+O(R^{-1})\\
\label{e2.38} & = & \sum_{q=1}^{9}M_{jq}\int\limits_{\mathbb{B}_R}D_{p(-)}(\xi)D_q(\nabla_\xi)^{\top}\delta(\xi)d\xi+O(R^{-1})\\
\nonumber & = & \sum_{q=1}^{9}M_{jq}D_q(-\nabla_\xi)D_{p(-)}(\xi)^{\top}|_{\xi=0}+O(R^{-1}) \\
&=&\left\{
\begin{array}{lll}
-M_{jp} & \mathrm{for} & p=1,\dots,6,\\
M_{jp} & \mathrm{for} & p=7,8,9
\end{array}
\right\}+O(R^{-1})=-(M_{(=)})_{jp}+O(R^{-1}).\nonumber
\end{eqnarray}
Here we have used that, first, $D_{p(-)}(\xi)$ is linear in $\xi$ and, therefore, $L_{(\top)}^{0}D_{p(-)}^\top=0$ and, second,
\begin{equation*}
L^0(\nabla_\xi)\Sigma^j(\xi):=\sum_{q=1}^{9}M_{jq}L^0(\nabla_\xi)(D_q(-\nabla_\xi)\Phi(\xi)^{\top})^{\top}=\sum_{q=1}^{9}M_{jq}D_q(-\nabla_\xi)^{\top}
\delta(\xi)
\end{equation*}
caused by the formula $L^0(\nabla_\xi)\Phi(\xi)=\delta(\xi)\mathbb{I}_{4\times4}$, i.e., by the definition of the fundamental matrix $\Phi$.

Let us process the left-hand  side $I_{left}^{jp}$ of
\eqref{e2.375}. Again integrating by parts, this time in the
domains $\Xi$ and $\omega$, it follows that
\begin{eqnarray}
I_{left}^{jp} & = & (N^0W^j,W_{(-)}^{p})_{\partial\Xi}-(N^0D_j^{\top},D_{p(-)}^\top)_{\partial\omega}\nonumber \\
& = & Q^0(W^j,W_{(-)}^{p};\Xi)+Q^0(D_j^{\top},D_{p(-)}^\top;\omega)\label{e2.39}\\
& = & Q_{(-)}^{0}(W^j,W^p;\Xi)+(A_{(-)}^{0})_{jp}mes_3\omega, \nonumber
\end{eqnarray}
where $mes_3\omega$ is the volume of $\omega$. Note that, first,
the equality $N^0W^j=-N^0D_{j}^\top$ on $\partial\omega$ is
inherited from \eqref{e2.31}   and \eqref{e2.13}, second,
$n^\omega$ and $-n^\omega$ imply the outward normals with respect
to the sets $\Xi$ and $\omega$, respectively, and, third,
\begin{eqnarray}
\nonumber Q^0(u,v_{(-)};\Xi) & = & (A^0D(\nabla_\xi)u,D(\nabla_\xi)v_{(-)})_{\Xi}\\
& = & (A_{(-)}^{0}D(\nabla_\xi)u,D(\nabla_\xi)v)_{\Xi}=Q_{(-)}^{0}(u,v;\Xi),\\
\nonumber Q_{(-)}^{0}(D_j^{\top},D_p^{\top};\omega) & = & (A_{(-)}^{0}\mathbf{e}_j,\mathbf{e}_p)_{\omega}=(A_{(-)}^{0})_{jp}mes_3\omega.
\end{eqnarray}

Comparing \eqref{e2.38}   and \eqref{e2.39}, we send $R$ to $+\infty$ and obtain the desired relation \eqref{e2.355}.$\blacksquare$

Theorem \ref{t2.8}  ensures the matrix $M_{(=)}$ in \eqref{e2.36}
to be symmetric, in particular,
$M^{\sf{ME}}=-(M^{\sf{EM}})^{\top}$. However, in contrast to the
polarization matrix in elasticity (cf. \cite{na100, na280,
NazSokSPN}) neither $M_{(=)}$, nor $M$ enjoy the
positivity/negativity property. In the case
$A^{\sf{ME}}=\mathbb{O}_{6\times3}$ the piezoelectricity problem
decouples into the elasticity and electricity problems so that,
\begin{equation}\label{e2.41}
M^{\sf{MM}}<0,\ M^{\sf{EE}}>0,\ M^{\sf{ME}}=-(M^{\sf{EM}})^{\top}=\mathbb{O}_{6\times3},
\end{equation}
provided, e.g., $mes_3\omega>0$. We emphasize that in
\eqref{e2.41}   $M^{\sf{EE}}$ is but the virtual mass tensor (see
\cite{PoSe}). By the perturbation argument, the matrix $M$ has six
negative and three positive eigenvalues, if the matrix
$A^{\sf{ME}}$ is sufficiently small (cf. Section
\ref{sec:examples}). However, for arbitrary $A^{\sf{ME}}$, this
property is still an open question.

We have examined the first asymptotic term \eqref{e2.32}   of  the
boundary layer type in the asymptotic ansatz \eqref{e2.17}.
By the representation \eqref{e2.34} (see Remark \ref{r2.5}), we
write the expansion of $w^1(\xi)$ for $\xi\rightarrow+\infty$ in
the matrix form as follows
\begin{equation}\label{e2.411}
w^1(\xi)=(D(\nabla_x)\Phi(\xi)^{\top})^{\top}M^{\top}\varepsilon^0+\widetilde{w}^1(\xi).
\end{equation}
The remainder $\widetilde{w}^1$ obeys the estimates \eqref{e2.288}.
\begin{remark}\label{r2.555}
Formula \eqref{e2.411}   can be derived in the following way:
\begin{eqnarray}
W^j(\xi)&=&\sum_{p=1}^{9}M_{jp}\sum_{q=1}^{3}\frac{\partial\Phi}{\partial\xi_q}(\xi)D_p(e_q)^\top+\widetilde{W}^j(\xi)\nonumber\\
&=&\sum_{p=1}^{9}M_{jp}\left(\sum_{q=1}^{3}D_p(e_q)\frac{\partial\Phi}{\partial\xi_q}(\xi)^\top\right)^\top+\widetilde{W}^j(\xi)\nonumber\\
&=&\left(\sum_{p=1}^{9}M_{jp}D_p(\nabla_\xi)\Phi(\xi)^\top\right)^\top+\widetilde{W}^j(\xi)\ .\nonumber
\end{eqnarray}
\end{remark}

\subsection{The second term in the boundary layer}
\label{4sec:2}
The system of differential equations in $\Xi$ in the exterior problem \eqref{e2.21} for the boundary layer $w^1$ is homogeneous. This leads to relatively simple formulae \eqref{e2.32} and \eqref{e2.411} for $w^1$. However, $w^2$ is determined from problem \eqref{e2.22} which enjoys the inhomegeneities both, in the boundary conditions and in the differential equations. Hence, the immediate objective becomes an inspection of the right-hand side $-L^\prime w^1$ for possible compensation and furthermore, an application of the same procedure as it is described in Section \ref{3sec:2}. However, the resulting asymptotic form \eqref{e2.50} of $w^2$ looks quite different compared to \eqref{e2.411}.

By virtue of \eqref{e2.14}   and \eqref{e2.15}, the operator
\begin{equation}\label{e2.42}
L^{\prime}(\xi,\nabla_\xi)=D(-\nabla_\xi)^{\top}\left(\sum_{j=1}^{3}\xi_jA^jD(\nabla_\xi)\right)
\end{equation}
gets the following homogeneity property:
\begin{equation}\label{e2.43}
L^{\prime}(\xi,\nabla_\xi)\rho^{\lambda}\varphi(\theta)=\rho^{\lambda-1}\psi(\theta),\ \xi\in\mathbb{R}^3\setminus\{0\}.
\end{equation}
Here $\lambda\in\mathbb{R}$, $(\rho,\theta)$ are the spherical coordinates in $\mathbb{R}^3$, $\rho=|\xi|$ and $\theta=\rho^{-1}\xi\in\mathbb{S}_1$,
and $\varphi,\psi\in C^{\infty}(\mathbb{S}_1)^4$.
Thus, by means of \eqref{e2.32}   and \eqref{e2.34}, \eqref{e2.288}, \eqref{e2.27}, we obtain that
\begin{equation}\label{e2.44}
F^{\prime}(\xi)=-L^{\prime}(\xi,\nabla_\xi)w^1(\xi)=D(\nabla_\xi)^{\top}(\rho^{-2}\Psi(\xi))+O(\rho^{-2}),\ \rho\rightarrow+\infty\ ,
\end{equation}
while the formula can be differentiated under the standard
convention $\nabla_xO(\rho^{-\lambda})=O(\rho^{-\lambda-1})$.
Due to the definition \eqref{e2.265}   of the Kondratiev norm the
right-hand side of \eqref{e2.44}   gives rise to the continuous
functional
\begin{equation*}
\begin{array}{c}
V_{0}^{1}(\Xi)^4\ni v\mapsto(F',v)_{\Xi},\\
|(F',v)_\Xi|\leq
c\displaystyle\int\limits_{\Xi}\rho^{-3}|v(\xi)|d\xi\leq
c\left(\int\limits_\Xi\rho^{-4}d\xi\right)^{1/2}\|\rho^{-1}v;L^{2}(\Xi)\|\leq
C\|v;V_0^1(\Xi)\|.
\end{array}
\end{equation*}
Thus, similarly to Proposition \ref{p2.3}, we obtain the existence
of a unique solution $w^2\in V_0^1(\Xi)^4$ to the problem \eqref{e2.22}. Now,
we  need to examine the behavior of $w^2(\xi)$ as
$\rho\rightarrow+\infty$. According to \cite{Ko} (see also
\cite[\S3.5]{NaPl}), first of all,  we have to determine the
power-law solution
\begin{equation}\label{e2.45}
Z(\xi)=\rho^{-1}\mathcal{Z}(\theta)
\end{equation}
to the system of differential equations
\begin{equation}\label{e2.46}
L^0(\nabla_\xi)Z(\xi)=\rho^{-3}\mathcal{F}(\theta):=D(\nabla_x)^{\top}(\rho^{-2}\Psi(\theta)),\
\xi\in\mathbb{R}^{3}\setminus\{0\},
\end{equation}
with the right-hand side taken from \eqref{e2.44}. Note that, in
general, the multiplier $\mathcal{Z}$ in \eqref{e2.45}   may be
linear in $\ln\rho$ but, owing to a special form of $\mathcal{T}$,
the next lemma proves the absence of the logarithm.
\begin{lemma}\label{l2.9}
The system \eqref{e2.46}   admits the power-law solution of form
\eqref{e2.45}, whose angular part $\mathcal{Z}(\theta)$ is defined
up to the linear combination
$c_1\Phi^1(\theta)+\dots+c_4\Phi^4(\theta)$, where
$c_j\in\mathbb{R}$ are arbitrary and $\Phi^j(\theta)$ is the trace on the unit
sphere $\mathbb{S}_1$ of the  column $\Phi^j(\xi)$ in the
fundamental matrix $\Phi$.
\end{lemma}

{\bf Proof.} After separation of variables and rewriting the
operator
$L^0(\nabla_\xi)=\rho^{-2}\mathfrak{L}(\theta,\nabla_\theta,\rho\partial_\rho)$
in the spherical coordinates $(\rho,\theta)$, the system
\eqref{e2.46} takes the form
\begin{equation}\label{e2.47}
\mathfrak{L}(\theta,\nabla_\theta,-1)\mathcal{Z}(\theta)=\mathcal{F}(\theta),\ \theta\in\mathbb{S}_1.
\end{equation}
By the Fredholm  alternative, this system on the unit sphere has a
solution if and only if the right-hand side $\mathcal{F}$ is
orthogonal to all solutions of the formally adjoint  homogeneous
system. Owing to \cite{MaPl1} (see also \cite[Lemma 3.5.9]{NaPl}),
the formally adjoint operator for
$\mathfrak{L}(\theta,\nabla_\theta,-1)$ is nothing but
$\mathfrak{L}_{(\top)}(\theta,\nabla_\theta,0)$, where
\begin{equation}\label{e2.477}
\rho^{-2}\mathfrak{L}_{(\top)}(\theta,\nabla_\theta,\rho\partial_\rho)=L_{(\top)}^{0}(\nabla_\xi)=L^0(\nabla_\xi)^{\ast}.
\end{equation}
By virtue of the polynomial property \eqref{e2.25}, any power-law
solution $X(\xi)=\rho^0\mathcal{X}(\xi)$ of
$L_{(\top)}^{0}(\nabla_\xi)X=0$ in $\mathbb{R}^3\setminus\{0\}$ is
a constant column in $\mathbb{R}^4$. Thus,  it suffices to verify
the orthogonality condition
\begin{equation}\label{e2.48}
\int\limits_{\mathbb{S}_1}\mathcal{F}(\theta)ds_\theta=0\in\mathbb{R}^4.
\end{equation}
Let $R>r>0$ and let $\Theta$ be the annulus $\{\xi:r<\rho<R\}$.
We have
\begin{gather*}
\displaystyle\ln\left(\frac{R}{r}\right)\int\limits_{\mathbb{S}_1}\mathcal{F}(\theta)ds_\theta=\int\limits_r^R\rho^{-1}d\rho\int\limits_{\mathbb{S}_1}\mathcal{F}(\theta)ds_\theta=\int\limits_\Theta\rho^{-3}\mathcal{F}(\theta)d\xi
=\displaystyle\int\limits_\Theta
D(\nabla_\xi)^{\top}(\rho^{-2}\psi(\theta))d\xi
\\
=\displaystyle\int\limits_{\mathbb{S}_R}D(\rho^{-1}\xi)^{\top}(\rho^{-2}\Psi(\theta))ds_\xi-\displaystyle\int\limits_{\mathbb{S}_r}D(\rho^{-1}\xi)^{\top}(\rho^{-2}\Psi(\theta))ds_\xi=0.
\end{gather*}
We have  used here the Gauss formula and the fact that the
integrands at $\rho=R$ and $\rho=r$ are equal to
$R^{-2}D(\theta)^{\top}\Psi(\theta)$ and
$r^{-2}D(\theta)^{\top}\Psi(\theta)$, respectively, so that the
integrals cancel each other.

Thus, the compatibility condition \eqref{e2.48}  holds true and
the system \eqref{e2.47}   admits a solution. It remains to recall
that any power-law solution \eqref{e2.45}   of the homogeneous
system \eqref{e2.46}   becomes a linear combination of the
fundamental matrix columns.$\blacksquare$

To assure the uniqueness of the solution \eqref{e2.45}, we impose
the condition
\begin{equation}\label{e2.49}
\int\limits_{\mathbb{S}_1}D(\theta)^{\top}A^0\mathfrak{D}(\theta,\nabla_\theta,-1)\mathcal{Z}(\theta)ds_\theta=0\in\mathbb{R}^4,
\end{equation}
where
$\rho^{-1}\mathfrak{D}(\theta,\nabla_\theta,\rho\partial_\rho)$ is
the matrix operator $D(\nabla_x)$ written, similarly to
\eqref{e2.477}, in the spherical coordinates $(\rho,\theta)$.

Now, we are in position to write an expansion at infinity for the
second boundary layer term in \eqref{e2.17}.
\begin{proposition}\label{p2.10}
The solution $w^2\in V_0^1(\Xi)^4$ of the problem \eqref{e2.22}
admits the asymptotic form
\begin{gather}\label{e2.50}
w^2(\xi)=Z(\xi)+\Phi(\xi)C+\widetilde{w}^2(\xi),\\
\label{XXX} |\nabla_\xi^k\widetilde{w}^2(\xi)|\leq
c_{k,\beta}\rho^{-2-k+\beta},\ k\in\mathbb{N}_0,\
\xi\in\mathbb{R}^3\setminus \mathcal{V},
\end{gather}
where $\beta>0$ is arbitrary, $Z$ is a power-law solution of form \eqref{e2.45}   and $C\in\mathbb{R}^4$ is determined as follows:
\begin{align}
C=-f(0)mes_3\omega+J\in\mathbb{R}^4,\label{e2.51}\\
J=\int\limits_{\mathbb{S}_1}D(\theta)^{\top}\sum_{j=1}^{3}\xi_jA^jD(\nabla_\xi)(D(\nabla_\xi)\Phi(\xi)^{\top})^{\top}ds_\xi M^\top\varepsilon^0.\label{e2.511}
\end{align}
\end{proposition}

{\bf Proof.} The asymptotic expansion \eqref{e2.50}   with a
certain column $C$ and the estimates \eqref{XXX}   result from
\cite{Ko} and \cite{MaPl2}, respectively (see also
\cite[Ch.3]{NaPl}). We again employ the method proposed in
\cite{MaPl1} to evaluate the constant column $C$.
Now, we use the Green formula in $\Xi_R$ for $w^2$ and $\sf{e}_p=(\delta_{p,1},\dots,\delta_{p,4})^{\top}$.
Recalling \eqref{e2.22}, we have
\begin{equation}\label{e2.52}
\begin{array}{l}
I_{left}:=-\displaystyle\int\limits_{\Xi_R}{\sf e}_p^{\top}L'w^1d\xi-\int\limits_{\partial\omega}{\sf e}_p^{\top}N'w^1ds\xi-\int\limits_{\partial\omega}{\sf e}_p^{\top}N'D(\xi)\varepsilon^0ds_\xi-\int\limits_{\partial\omega}{\sf e}_p^{\top}N'Uds_\xi\\
=\displaystyle\int\limits_{\Xi_R}{\sf e}_p^{\top}L^0w^2d\xi+\int\limits_{\partial\omega}{\sf e}_p^{\top}N^0w^2ds\xi=\int\limits_{\mathbb{S}_R}{\sf e}_p^{\top}N^0w^2ds\xi=:I_{right}.
\end{array}
\end{equation}
Here $N^0(\xi,\nabla_\xi)=D(\theta)^{\top}A^0D(\nabla_\xi)$ on the sphere $\mathbb{S}_R$ with the unit normal vector $\theta=\rho^{-1}\xi$ (cf. \eqref{e2.20}   and \eqref{e2.49}).
Similarly to the calculation \eqref{e2.38},  using \eqref{e2.50}   and \eqref{e2.49}, we get
\begin{equation}\label{e2.53}
\begin{array}{c}
I_{right}=-\int\limits_{\mathbb{S}_R}{\sf e}_p^{\top}N^0\mathcal{Z}ds_\xi-\int\limits_{\mathbb{S}_R}{\sf e}_p^{\top}N^0\Phi ds_\xi C+O(R^{-1})=\\
=\displaystyle\int\limits_{\mathbb{B}_R}{\sf e}_p^{\top}L^0\Phi d\xi C+O(R^{-1})=C_p+O(R^{-1}).
\end{array}
\end{equation}
By integrating by parts, the last couple of integrals in $I_{left}$ turns into
\begin{equation}\label{e2.54}
-\displaystyle\int\limits_{\Xi_R}{\sf e}_p^{\top}L'w^1d\xi-\int\limits_{\partial\omega}{\sf e}_p^{\top}N'w^1ds\xi=\int\limits_{\mathbb{S}_R}{\sf e}_p^{\top}D(\theta)^{\top}\sum_{j=1}^{3}\xi_jA^jD(\nabla_\xi)(D(\nabla_\xi)\Phi(\xi)^{\top})^{\top}ds_\xi M\varepsilon^0+O(R^{-1}).
\end{equation}
Here we have applied the decomposition \eqref{e2.411}   of $w^1$
together with the estimate \eqref{e2.288}   for the remainder.
Since its integrand is a positive homogeneous function in $\xi$ of
degree $-2$ (cf. \eqref{e2.27}) the integral $J_p$ over
$\mathbb{S}_R$ in \eqref{e2.54}   is independent of the radius $R$
and becomes an entry of column \eqref{e2.511}.

The first couple of integrals in \eqref{e2.52}   is equal to
\begin{equation*}
\begin{array}{l}
\displaystyle\int\limits_{\partial\omega}{\sf e}_p^{\top}N'D(\xi)\varepsilon^0ds_\xi-\int\limits_{\partial\omega}{\sf e}_p^{\top}N^0Uds_\xi=-\int\limits_{\omega}{\sf e}_p^{\top}(L'D(\xi)\varepsilon^0+L^0U)d\xi\\
=-mes_3\omega{\sf e}_p^{\top}(-\sum_{j=1}^{3}D(e_j)^{\top}A^j\varepsilon^0+L^0(\nabla_\xi)U(\xi))=-f_p(0)mes_3\omega.
\end{array}
\end{equation*}
Here, the elementary formula \eqref{e2.16}   has been taken into account.

Now the limit passage $R\rightarrow+\infty$ in \eqref{e2.52}-\eqref{e2.54}   furnishes \eqref{e2.51}   and \eqref{e2.511}.$\blacksquare$
\begin{remark}\label{r2.0}
Proposition \ref{p2.6}  can be proved by an application of the
method \cite{MaPl1} in the same way as it is made in Proposition
\ref{p2.10} and Theorem \ref{t2.8}. We only mention that the
columns $d^1,\dots,d^7$ of the matrix $d(\xi)$ in \eqref{e2.12}
satisfy simultaneously the homogeneous problem \eqref{e2.21} and
the formally adjoint boundary value problem in $\Xi$ with the
operators $L_{(\top)}^0(\nabla_\xi)$ and
$N_{(\top)}^0(\xi,\nabla_\xi)$, respectively.
\end{remark}

\subsection{The regular correction term}
\label{5sec:2}

Let us consider now the subsequent term in the asymptotic ansatz
\eqref{e2.17}, namely the regular correction term $\mathbf{u}(x)$.

By means of \eqref{e2.411}   and \eqref{e2.50}, we have
\begin{equation}\label{e2.55}
\begin{array}{l}
hw^1(h^{-1}x)+h^2w^2(h^{-1}x)=h(S^2(h^{-1}x)+\widetilde{w}^1(h^{-1}x))+\\
+h^2(S^1(h^{-1}x)+\widetilde{w}^2(h^{-1}x))=h^3(S^2(x)+S^1(x))+O(h^4(|x|^{-3}+|x|^{-2}))
\end{array}
\end{equation}
where, according to \eqref{e2.27}   and \eqref{e2.45}, we have set
\begin{gather}\label{e2.56}
S^2(\xi)=(D(-\nabla_\xi)\Phi(\xi)^{\top})^{\top}M^{\top}\varepsilon^0,\ S^1(\xi)=Z(\xi)+\Phi(\xi)C,\\
\notag S^p(t\xi)=t^{-p}S^p(\xi).
\end{gather}
Therefore, this is $h^3\mathbf{u}(x)$ in the asymptotic ansatz
\eqref{e2.17}   that compensates the main part of a discrepancy
produced by the boundary layer terms $w^1$ and $w^2$.

Taking into account the equalities  $L^0S^2=0$ and $L^0S^1=-L'S^2$
designated in two last sections, we arrive at the following
representation of the discrepancy in the system \eqref{e2.2}  :
\begin{equation}\label{e2.57}
\begin{array}{l}
\mathbf{f}(x)=-L(x,\nabla_x)(\chi(x)(S^2(x)+S^1(x)))=\\
=-[L,\chi](S^2(x)+S^1(x))-\chi(x)(L(x,\nabla_x)-L^0(\nabla_x)-L'(x,\nabla_x))S^2(x)-\\
-\chi(x)(L(x,\nabla_x)-L^0(\nabla_x))S^1(x).
\end{array}
\end{equation}
Here $[L,\chi]$ stands for the commutator of the differential operator $L$ and the cut-of function $\chi$, i.e.,
\begin{equation}\label{e2.99}
[L,\chi]=D(-\nabla_x)^{\top}A(x)D(\nabla_x\chi(x))-D(\nabla_x\chi(x))^{\top}A(x)D(\nabla_x).
\end{equation}
Recalling \eqref{e2.14}   and \eqref{e2.15}, in view of \eqref{e2.56}, we obtain
that
\begin{equation}\label{e2.58}
|\mathbf{f}(x)|\leq c|x|^{-2}.
\end{equation}

We see that the regular correction term $\mathbf{u}$ must satisfy the piezoelectricity problem
\begin{align}
D(-\nabla_x)^{\top}A(x)D(\nabla_x)\mathbf{u}(x)=\mathbf{f}(x),\ x\in\Omega,\label{e2.59}\\
D(n(x))^{\top}A(x)D(\nabla_x)\mathbf{u}(x)=0,\ x\in\Gamma_\sigma,\ \mathbf{u}(x)=0,\ x\in\Gamma_u.\label{e2.60}
\end{align}
We emphasize that the sum
$h\widetilde{w}^1(h^{-1}x)+h^2\widetilde{w}^2(h^{-1}x)$ in
\eqref{e2.55}   becomes of order $h^4$ only at a distance from the
coordinate origin $x=0$. However, we have extended  equations
\eqref{e2.59}   over the whole domain $\Omega$ because the
singularity $O(|x|^{-2})$ of the right-hand side $\mathbf{f}(x)$
is not too strong. In particular, by \eqref{e2.58}, the functional
on the right-hand side  in the integral identity
\begin{equation}\label{e2.61}
Q(\mathbf{u},\mathbf{v};\Omega)=(\mathbf{f},\mathbf{v})_\Omega,\ \mathbf{v}\in\mathring{H}^1(\Omega;\Gamma_u)^4,
\end{equation}
serving for the problem \eqref{e2.59}, \eqref{e2.60}   (cf. \eqref{e1.10}), is continuous due to the estimate
\begin{equation*}
\begin{array}{l}
|(\mathbf{f},\mathbf{v})_\Omega|\leq c\left(\displaystyle\int\limits_{\Omega}|x|^2|\mathbf{f}(x)|^2dx\right)^{1/2}\left(\displaystyle\int\limits_\Omega|x|^{-2}|\mathbf{v}(x)|^2dx\right)^{1/2}\leq \\
\leq c\left(\displaystyle\int\limits_0^{diam\Omega}r^2r^{-4}r^2dr\right)^{1/2}\|\nabla_x\mathbf{v};L_2(\Omega)\|\leq C\|\mathbf{v};H^1(\Omega)\|
\end{array}
\end{equation*}
and the one-dimensional Hardy inequality mentioned above
\eqref{e2.265}. Hence, in the analogy with Proposition \ref{p1.1},
the Lax-Milgram lemma ensures the existence and uniqueness of the
solution $\mathbf{u}\in\mathring{H}^1(\Omega;\Gamma_u)^4$. These
observations complete the evaluation of all asymptotic terms
detached in \eqref{e2.17}.
\begin{remark}\label{r2.12}
The singularity of $\mathbf{f}$ can lead to a logarithmical singularity of the solution $\mathbf{u}$.
However, we shall need only the following inequalities with arbitrary $\beta>0$\ :
\begin{equation}\label{e2.68}
|\mathbf{u}(x)|\leq c_{\beta}|x|^{-\beta},\ |\nabla_x\mathbf{u}(x)|\leq c_{\beta}|x|^{-1-\beta}
\end{equation}
delivered by a result in \cite{MaPl2} (see also \cite[\S
3.6]{NaPl}).
\end{remark}

For the further usage, it is convenient to rewrite the  ansatz
\eqref{e2.17}   in a different form, namely
\begin{equation}\label{e2.62}
u^h(x)=u(x)+h^3\mathbf{U}(x)+\chi(x)(h\widetilde{w}^1(h^{-1}x)+h^2\widetilde{w}^2(h^{-1}x))+\widetilde{u}^h(x),
\end{equation}
where, in accordance with \eqref{e2.55}   and \eqref{e2.56},
\begin{equation}\label{e2.63}
\mathbf{U}(x)=\mathbf{u}(x)+\chi(x)(S^2(x)+S^1(x)).
\end{equation}
In other words, we detach $hS^2(h^{-1}x)$ and $h^2S^1(h^{-1}x)$ from the boundary layer terms and attach them to the regular term $\mathbf{u}$.
Therefore, the remainder $\widetilde{u}^h$ in \eqref{e2.62}   stays the same as in the original ansatz \eqref{e2.17}.

Let us derive an {\it almost explicit} formula for \eqref{e2.63}. To
this end, let $G(x,y)$ be the Green matrix for the
piezoelectricity problem \eqref{e1.7}-\eqref{e1.9}, i.e.,
\begin{equation}\label{e2.64}
\begin{array}{l}
D(-\nabla_x)^{\top}A(x)D(\nabla_x)G(x,y)=\delta(x-y)\mathbb{I}_{4\times4},\ x\in\Omega, \\
D(n(x))^{\top}A(x)D(\nabla_x)G(x,y)=0,\ x\in\Gamma_\sigma,\ u(x)=0,\ x\in\Gamma_u
\end{array}
\end{equation}
Of course, the relations \eqref{e2.64} are understood  in the  sense
of distributions, so that, $G\in L^2(\Omega)^{4\times4}$, $G\in
L^2(\partial\Omega)^{4\times4}$ and
\begin{equation*}
(G,L_{(\top)}v)_\Omega+(G,N_{(\top)}v)_{\Gamma_\sigma}=v(y),\ v\in C_{c}^{\infty}(\overline{\Omega};\Gamma_u)^4,
\end{equation*}
where the linear space
$C_{0}^{\infty}(\overline{\Omega};\Gamma_u)$ consists of
infinitely differentiable functions in $\overline{\Omega}$ which
vanish on $\Gamma_u$. Since $A$ is a smooth matrix function inside
of the ball $\mathbb{B}_R$ (see \eqref{e2.6}), the Green matrix is
properly defined for $y\in\mathbb{B}_R$ (see \cite{GeShi, Herm})
and
\begin{equation*}
(x\mapsto G(x,y)-\Phi(x,y))\in H^1(\Omega)^{4\times4}.
\end{equation*}
Moreover, $G$ can be differentiated in the second argument and we set
\begin{equation}\label{e2.65}
G^0(x)=G(x,0),\ \mathbf{G}^0(x)=D(-\nabla_y)G(x,y)|_{y=0}.
\end{equation}
By repeating the considerations in and around of Lemma \ref{l2.9}, we can detect that
\begin{equation}\label{e2.66}
G^0-\Phi\in H^1(\Omega)^{4\times4},\ \mathbf{G}^0-D(\nabla_x)\Phi-\mathbf{Z}-\mathbf{K}\Phi\in H^1(\Omega)^{9\times4},
\end{equation}
where $\mathbf{K}$ is a certain  matrix  of the
size $9\times4$ with real entries and $\mathbf{Z}$ is such that
$Z(x)=\mathbf{Z}(x)M^{\top}\varepsilon^0$ (cf. \eqref{e2.42}   and
\eqref{e2.44}-\eqref{e2.46}). Since, by definition of $\mathbf{u}$
and $S^q$, the vector function $\mathbf{U}$ verifies the boundary
conditions \eqref{e2.60}   and the homogeneous system
\eqref{e2.59} everywhere in $\Omega$, except at the point
$\mathcal{O}$. Let us now compare singularities in \eqref{e2.66}
and \eqref{e2.63} to conclude that
\begin{equation}\label{e2.67}
\mathbf{U}(x)=\mathbf{G}^0(x)M^{\top}\varepsilon^0-G^0(x)f(0)mes_3\omega.
\end{equation}

\begin{remark}\label{r2.11}
1.\ We emphasize that the differential operator $D(-\nabla_y)$ in \eqref{e2.65}   is replaced by $D(\nabla_x)$ in \eqref{e2.66}.
This is due to the evident relationship $D(-\nabla_y)\Phi(x-y)=D(\nabla_x)\Phi(x-y)$.\\

2.\ If $A$ is a constant matrix, then  the terms $\mathbf{Z}$ and
$\mathbf{K}\Phi$ are absent in \eqref{e2.66}, in other words,
their presence results from the variable  coefficients of
differential operator \eqref{e2.15}. Therefore, the column
$-f(0)mes_3\omega$ occurs on the right-hand side of \eqref{e2.67}
instead of the column \eqref{e2.51}. To ensure that the additional
column \eqref{e2.511}   does not effect the form of  the last term
in \eqref{e2.67}, one may put $\varepsilon^0=0$ to see that then
$J=0$. A direct calculation leading to formula \eqref{e2.67}   can
be found in \cite{NaSoElast} for the three-dimensional elasticity
problem.
\end{remark}

Since the coordinate origin $\mathcal{O}$ is situated inside
$\omega_h$, i.e., outside $\overline{\Omega}_h$ (cf. Section
\ref{1sec:2}), the second term \eqref{e2.67}   in the new ansatz
\eqref{e2.62}   is smooth in the domain $\Omega(h)$, although the
Green matrices \eqref{e2.65}   have singularities at
$\mathcal{O}$.

\section{Justification of asymptotics and analysis of shape functionals}
\label{sec:3}

\subsection{The justification of asymptotics}
\label{1sec:3}

The difference
\begin{equation}\label{e3.1}
\widetilde{u}^h=u^h-u-\chi(hw^1+h^2w^2)-h^3\mathbf{u}
\end{equation}
(see \eqref{e2.17}   and \eqref{e2.62}) satisfies the integral identity
\begin{equation}\label{e3.2}
Q(\widetilde{u}^h,v;\Omega(h))=\widetilde{\mathcal{F}}^h(v),\ v\in\mathring{H}^1(\Omega(h);\Gamma_u)^4,
\end{equation}
where $\widetilde{\mathcal{F}}^h$ is a certain functional (see, e.g., \eqref{e3.5}).
If the estimate
\begin{equation}\label{e3.3}
|\widetilde{\mathcal{F}}^h v|\leq ch^{\alpha+5/2}\|v;H^1(\Omega(h))\|
\end{equation}
is proved, we could take $v=\widetilde{u}^h$ in order to conclude by using \eqref{e1.000}   that
\begin{equation}\label{e3.4}
\|\widetilde{u}^h;H^1(\Omega(h))\|\leq ch^{\alpha+5/2}.
\end{equation}
In the sequel, it is shown, that the constants in \eqref{e3.3} and \eqref{e3.4} are independent of the small parameter $h$.

To verify \eqref{e3.3}, first, we assume that $v$ vanishes in the ball $\mathbb{B}_{2R/3}$, therefore,  $\chi v=0$.
Then, we have
\begin{equation}\label{e3.5}
\widetilde{\mathcal{F}}^h(v)=Q(u^h-u-h^3\mathbf{u},v;\Omega(h))=Q(u^h,v;\Omega(h))-Q(u,v;\Omega)-h^3Q(\mathbf{u},v;\Omega).
\end{equation}
Recalling  \eqref{e2.555}, \eqref{e1.10}   and \eqref{e2.61}, we
observe that the support  of the vector function \eqref{e2.57}
satisfies $supp\,\mathbf{f}\subset\overline{\mathbb{B}}_{2R/3}$
(each term in \eqref{e2.57} contains either a cut-off function $\chi$ supported in the ball, or its derivatives)
and, hence, \eqref{e3.5}   is null.

Second, let $supp\,v\subset\mathbb{B}_R\setminus\omega_h$.
We write
\begin{equation}\label{e3.6}
\begin{array}{l}
\widetilde{F}^h(v)=(f,v)_{\Omega(h)}-(AD(\nabla_x)u,D(\nabla_x)v)_{\Omega(h)}-h^3(AD(\nabla_x)\mathbf{u},D(\nabla_x)v)_{\Omega(h)}\\
-h(AD(\nabla_x)(\chi w^1),D(\nabla_x)v)_{\Omega(h)}
-h^2(AD(\nabla_x)(\chi w^2),D(\nabla_x)v)_{\Omega(h)}\\
=:(f,v)_{\Omega(h)}-I^{u}-h^3I^{\mathbf{u}}-hI_1^w-h^2I_2^w.
\end{array}
\end{equation}

Since the vector functions $u$ and $\mathbf{u}$ are smooth in $\mathbb{B}_R\setminus\omega_h$, we integrate by parts and obtain
\begin{align}
I^{u}=(f,v)_{\Omega(h)}+(D(n^h)^{\top}A^0\varepsilon^0,v)_{\partial\omega_h}+(D(n^h)^{\top}A^0D(\nabla_x)U,v)_{\partial\omega_h}\nonumber\\ \label{e3.7}
+\sum_{j=1}^{3}(D(n^h)^{\top}x_jA^j\varepsilon^0,v)_{\partial\omega_h}+\widetilde{I}^{u},\\
\widetilde{I}^{u}=(D(n^h)^{\top}AD(\nabla_x)(u-D(x)^{\top}\varepsilon^0-U),v)_{\partial\omega_h}+\nonumber\\
(D(n^h)^{\top}(A-A^0-\sum_{j=1}^{3}x_jA^j)\varepsilon^0,v)_{\partial\omega_h}+(D(n^h)^{\top}(A-A^0)D(\nabla_x)U,v)_{\partial\omega_h},\nonumber\\
I^{\mathbf{u}}=(\mathbf{f},v)_{\Omega(h)}+\widetilde{I}^{\mathbf{u}},\quad \widetilde{I}^{\mathbf{u}}=(D(n^h)^{\top}AD(\nabla_x)\mathbf{u},v)_{\partial\omega_h}.\label{e3.8}
\end{align}

To process the terms $\widetilde{I}^{u}$ and $\widetilde{I}^\mathbf{u}$, we recall the inequality
\begin{equation}\label{e3.9}
\int\limits_{\Omega(h)}|x|^{-2}|v(x)|^2dx\leq c\|v;H^1(\Omega(h))\|^2\ ,
\end{equation}
which is a consequence of the one-dimensional Hardy inequality
(cf. \cite[\S4.5]{NaPl}) and the trace inequality (see \cite{Lad})
\begin{equation}\label{e3.10}
\int\limits_{\partial\omega_h}|v(x)|^2ds_x\leq ch\|v;H^1(\Omega(h))\|^2\ ,
\end{equation}
where the constants $c$ are independent of $h\in(0,h_0]$ and $v$.

Now by \eqref{e3.10}   and \eqref{e2.68}, we readily derive that
\begin{equation}\label{e3.11}
\begin{array}{rcl}
h^3|\widetilde{I}^{\mathbf{u}}|\leq ch^3h^{-1-\beta}\int\limits_{\partial\omega_h}|v(x)|ds_x&\leq& ch^{2-\beta}(mes_2\partial\omega_h)^{1/2}h^{1/2}\|v;H^1(\Omega(h))\|\\
&=&Ch^{-\beta+7/2}\|v;H^1(\Omega(h))\|.
\end{array}
\end{equation}
Analogously, by means of \eqref{e2.14}, \eqref{e2.10}   and \eqref{e3.10}, we have
\begin{equation}\label{e3.12}
|\widetilde{I}^{u}|\leq c(h^{1+\alpha}+h^2+h^2)\int\limits_{\partial\omega_h}|v(x)|ds_x\leq ch^{\alpha+5/2}\|v;H^1(\Omega(h))\|.
\end{equation}
We may choose $\beta=1-\alpha>0$ in order to equalize the final exponents of $h$ in \eqref{e3.11}   and \eqref{e3.12}.

Dealing with $I_2^w$, we write
\begin{equation}\label{e3.13}
\begin{array}{rcl}
I_2^w & = & (AD(\nabla_x\chi)S^1,D(\nabla_x)v)_{\Omega(h)}-(AD(\nabla_x)S^1,D(\nabla_x\chi)v)_{\Omega(h)}\\
&+&(A^0D(\nabla_x)w^2,D(\nabla_x)(\chi v))_{\Omega(h)}+\widetilde{I}_2^w,
\end{array}
\end{equation}
\begin{equation}\label{e3.14}
\begin{array}{rcl}
\widetilde{I}_2^w&=&(AD(\nabla_x\chi)(w^2-S^1),D(\nabla_x)v)_{\Omega(h)}-(AD(\nabla_x\chi)(w^2-S^1),D(\nabla_x\chi)v)_{\Omega(h)}\\
&+&((A-A^0)D(\nabla_x)w^2,D(\nabla_x)(\chi v))_{\Omega(h)}.
\end{array}
\end{equation}
Here, we detach  $S^1(h^{-1}x)$ from $w^2(h^{-1}x)$ (cf.
\eqref{e2.55}) and commute twice the differential operator
$D(\nabla_x)$ with the cut-off function $\chi$ (see
\eqref{e2.99}).

In view of \eqref{e2.14}   and \eqref{e2.50}, the absolute value
of the last expression in \eqref{e3.14}, multiplied by $h^2$
according to the definition of $I_2^w$ in \eqref{e3.6}, does not
exceed the sum of the following two expressions:
\begin{equation}\label{e3.15}
\begin{array}{l}
ch^2\displaystyle\int\limits_{\Omega\setminus\mathbb{B}_Rh}|x|h^{-1}\left(\displaystyle\frac{|x|}{h}\right)^{-3+\beta}|D(\nabla_x)(\chi(x)v(x))|dx\leq \\
\leq ch^{4-\beta}\left(\displaystyle\int\limits_{Rh}^{diam\Omega}r^2r^{-6+2\beta}r^2dr\right)^{1/2}\|v;H^1(\Omega(h))\|\leq ch^{7/2}\|v;H^1(\Omega(h))\|
\end{array}
\end{equation}
and
\begin{equation}\label{e3.155}
\begin{array}{l}
ch^2\displaystyle\int\limits_{\mathbb{B}_{Rh}\setminus\omega_h}|x||D(\nabla_x)\widetilde{w}^2(h^{-1}x)||D(\nabla_x)(\chi v)|dx\\
\leq ch^2Rh\left(\displaystyle\int\limits_{\mathbb{B}_R\setminus\omega}h^{-2}|D(\nabla_\xi)\widetilde{w}^2(\xi)|^2d\xi h^3\right)^{1/2}\|v;H^1(\Omega(h))\|\\
\leq ch^{7/2}\|v;H^1(\Omega(h))\|.
\end{array}
\end{equation}
The radius $R$ is chosen such that
$\mathbb{B}_R\supset\overline{\omega}$. Since the support of
$|D(\nabla_x\chi)|$ belongs to the annulus
$\overline{\mathbb{B}_{2R/3}}\setminus\mathbb{B}_{R/3}$ where,
according to \eqref{e2.50},
\begin{equation*}
|w^2(h^{-1}x)-S^1(h^{-1}x)|+|\nabla_x(w^2(h^{-1}x)-S^1(h^{-1}x))|\leq ch^{2-\beta},
\end{equation*}
the remaining  terms in \eqref{e3.14}, again after multiplication
by $h^2$, are bounded by $ch^{4-\beta}\|v;H^1(\Omega)\|$ while we
may set $\beta=1/2$ to achieve the same exponent as in
\eqref{e3.15}. In other words, for $p=2$, we now have
\begin{equation}\label{e3.16}
h^p|\widetilde{I}_p^w|\leq ch^{7/2}\|v;H^1(\Omega(h))\|.
\end{equation}
By formulae \eqref{e2.411}, \eqref{e2.288}   and \eqref{e2.55},
the similar argument  leads to the estimate \eqref{e3.16}   for
the remainder in the representation
\begin{equation}\label{e3.17}
\begin{array}{l}
I_1^w=(AD(\nabla_x\chi)S^2,D(\nabla_x)v)_{\Omega(h)}-(AD(\nabla_x)S^2,D(\nabla_x\chi)v)_{\Omega(h)}+\\
(A^0D(\nabla_x)w^1,D(\nabla_x)(\chi v))_{\Omega(h)}+\sum_{j=1}^{3}(x_jA^jD(\nabla_x)w^1,D(\nabla_x)(\chi v))_{\Omega(h)}+\widetilde{I}_1^w.
\end{array}
\end{equation}

Now, we are in position to conclude the estimate \eqref{e3.3}
for the functional $\widetilde{F}^h$ in \eqref{e3.2}, \eqref{e3.5}
and \eqref{e3.6}. To this end, we list several facts. First, the
inner product $(f,v)_{\Omega(h)}$ on the right hand-side of
\eqref{e3.6} cancels the same product in \eqref{e3.7}. Second, the
equality
\begin{equation*}
(\mathbf{f,v})_{\Omega(h)}=-(AD(\nabla_x\chi)(S^2-S^1),D(\nabla_x)v)_{\Omega(h)}+(AD(\nabla_x)(S^2-S^1),D(\nabla_x\chi)v)_{\Omega(h)}
\end{equation*}
is inherited from the definitions \eqref{e2.57}   and
\eqref{e2.99}. Third, we make the coordinate dilation
$x\mapsto\xi=h^{-1}x$ in the first couples of terms on the right
hand-side of \eqref{e3.13} and \eqref{e3.17}, simultaneously
multiplying the terms by $h^2$ and $h$, respectively. Noting that
$S^p(h^{-1}x)=h^pS^p(x)$, $p=1,2$, we see that these couples and
$h^3(\mathbf{f,v})_{\Omega(h)}$ annihilate. Finally, we recall
the integral identities \eqref{e2.26}, serving for the problems
\eqref{e2.21}   and \eqref{e2.22}, and after the substitutions
$x\mapsto\xi$ and $v(\xi)\mapsto\chi(h\xi)v(h\xi)$, we detect all
terms in the identities on the right hand-sides of \eqref{e3.7},
\eqref{e3.13}   and \eqref{e3.17}. Thus,
\begin{equation*}
\widetilde{F}^h(v)=\widetilde{I}^{u}+h^3\widetilde{I}^{\mathbf{u}}+h\widetilde{I}_1^w+h^2\widetilde{I}_2^w
\end{equation*}
and the inequality \eqref{e3.3}   holds true by virtue of \eqref{e3.12}, \eqref{e3.11}   and \eqref{e3.16}   with $p=1,2$.
We notice that the lowest exponent
$\alpha+5/2$ of $h$ occurs in \eqref{e3.12}   because $\alpha\in(1/2,1)$ and $\alpha+5/2\in(3,7/2)$.

We now formulate the result.
\begin{theorem}\label{t3.1}
Let all assumptions in Section \ref{1sec:2} be valid, in
particular, the inclusion \eqref{e2.8}   with $\alpha\in(1/2,1)$.
Then the solution $u^h$ of the piezoelectricity problem
\eqref{e2.2}-\eqref{e2.5}   and its approximation constructed in
Section \ref{sec:2} are in the relationship
\begin{equation}\label{e3.18}
\|u^h-u-h^3\mathbf{u}-\chi(hw^1+h^2w^2);H^1(\Omega(h))\|\leq ch^{\alpha+5/2}N,
\end{equation}
where the constant $c$ is independent of the parameter $h\in(0,h_0]$ and the right-hand sides $f$, $g$ while
\begin{equation}\label{e3.19}
N=\|f;L^2(\Omega)\|+\|g;L^2(\partial\Omega)\|+\|f;C^{2,\alpha}(\mathbb{B}_R)\|.
\end{equation}
In asymptotic estimate \eqref{e3.18} $u$ stands for a solution of piezoelectricity mixed boundary value problem \eqref{e1.7}-\eqref{e1.9} in the intact body $\Omega$ and $\mathbf{u}$ for the main regular corrector which is a solution of problem \eqref{e2.59}, \eqref{e2.60} in $\Omega$ and admits the representation \eqref{e2.63} with the smooth addendum $\mathbf{U}$ and the singular components \eqref{e2.56}. The boundary layer terms $w^1$ and $w^2$, which are as well present in   \eqref{e3.18}, are given by solutions to the exterior problems
\eqref{e2.21} and \eqref{e2.22},
decay at infinity and take asymptotic forms  \eqref{e2.411} and \eqref{e2.50}, respectively.
\end{theorem}
\begin{remark}\label{r3.2}
The obtained estimate \eqref{e3.18}   is asymptotically sharp, in
particular, it satisfies the "first omitted term"\, rule. Indeed,
for the smooth data $A$ and $f$, the subsequent asymptotic term in the
ansatz \eqref{e2.17}   is $h^3\chi(x)w^3(h^{-1}x)$, the
$H^1(\Omega(h))$-norm of the latter term is just $O(h^{7/2})$.
This bound appears in \eqref{e3.19}   if $\alpha\rightarrow1-0$.
Moreover, the estimate \eqref{e3.18}   holds true when the last
addendum in \eqref{e3.19}   is changed for
$\|f;C^{3,\alpha_1}(\mathbb{B}_R)\|$ with any $\alpha_1\in(0,1)$.
If the right-hand side $f\in C^{2,\alpha}(\mathbb{B}_R)^4$ in the
equations \eqref{e2.2}   is not sufficiently smooth, e.g.,
\begin{equation*}
f(x)=f^0(x)+|x|^{2+\alpha}f^1(\theta),\ f^0\in C^{\infty}(\mathbb{B}_R)^4,\ f^1\in C^{\infty}(\mathbb{S}_1)^4,
\end{equation*}
then the asymptotic ansatz \eqref{e2.17}   gains the boundary
layer term $h^{2+\alpha}\chi(x)w^{2+\alpha}(h^{-1}x)$ with the
Sobolev norm in $\Omega(h)$ of the some order $h^{\alpha+5/2}$ as
on the right hand-side of \eqref{e3.18}.
\end{remark}

A direct calculation show that
\begin{equation}
h^{j}\|\chi w^{j};H^{1}(\Omega(h))\|=O(h^{j+1/2}),\ j=1,2,
\end{equation}
and, therefore, in view of the relation $\alpha+5/2>3$ (see
\eqref{e2.8}), the $H^{1}(\Omega(h))$-norm of each of the detached
asymptotic terms in \eqref{e3.18}   (cf. \eqref{e2.17}   and
\eqref{e2.62}) is of order $h^{s}$ with $s\leq3$. In other
words, Theorem \ref{t3.1} justifies the constructed asymptotics of
solution $u^{h}$, indeed.

\subsection{The energy and the electric enthalpy}
\label{2sec:3}

We proceed with  energy functional \eqref{e1.13}, assuming for simplicity that the volume forces and the volume charges are absent, i.e., $f=0$ on the right hand-sides of \eqref{e1.7}   and \eqref{e2.2}.
Then, integrating by parts and taking   into account formulae \eqref{e2.62}   and \eqref{e3.4}, we have
\begin{equation}\label{e3.20}
\begin{array}{rcl}
\mathcal{U}(u^h;\Omega(h))&=&\frac{1}{2}(D(n)^{\top}AD(\nabla_x)u^h,u^h)_{\Gamma_\sigma}-(g,u^h)_{\Gamma_\sigma}=-\frac{1}{2}(g,u^h)_{\Gamma_\sigma} \\
&=&-\frac{1}{2}(g,u)_{\Gamma_\sigma}-\frac{1}{2}h^3(g,\mathbf{u})_{\Gamma_\sigma}+O(h^{\alpha+5/2})\ .
\end{array}
\end{equation}
Let $\mathfrak{v}^{\sf M}\in\mathring{H}^1(\Omega;\Gamma_u)^4$ and
$\mathfrak{v}^{\sf E}\in\mathring{H}^1(\Omega;\Gamma_u)^4$ imply
the solutions of the problem \eqref{e1.7}-\eqref{e1.9}   with the
right-hand sides
\begin{equation}\label{e3.21}
\mathfrak{g}^{\sf M}=(g_1^{\sf M},g_2^{\sf M},g_2^{\sf M},0)^{\top},\ \mathfrak{g}^{\sf E}=(0,0,0,g_4^{\sf E})^{\top}.
\end{equation}
Using the representation \eqref{e2.67}   with $f(0)=0$ and the modified column $\mathbf{U}_{(-)}$ (see \eqref{e2.365}), we obtain
\begin{equation}\label{e3.22}
\begin{array}{rcl}
(\mathfrak{g}^{\sf M},\mathbf{U})_{\Gamma_{\sigma}}&=&(\mathfrak{g}^{\sf M},\mathbf{U}_{(-)})_{\Gamma_\sigma}=(D(n)^{\top}AD(\nabla_x)\mathfrak{v}^{\sf M},\mathbf{U}_{(-)})_{\Gamma_\sigma}\\
&=&(\mathfrak{v}^{\sf M},D(-\nabla_x)^{\top}A^{\top}D(\nabla_x)\mathbf{U}_{(-)})_{\Gamma_\sigma}\\
&=&(\varepsilon^0)^{\top}M(\mathfrak{v}^{\sf M},D(-\nabla_x)^{\top}A^{\top}D(\nabla_x)\mathbf{G}_{(-)}^0)_{\Omega}\\
&=&(\varepsilon^0)^{\top}M(\mathfrak{v}^{\sf M},(D(\nabla_x)^\top\delta(x))_{(-)})_\Omega=-(\varepsilon^0)^{\top}M\mathfrak{e}_{(-)}^{\sf M},\\
(\mathfrak{g}^{\sf E},\mathbf{U})_{\Gamma_{\sigma}}&=&(\mathfrak{g}^{\sf E},\mathbf{U}_{(-)})_{\Gamma_\sigma}=-(\mathfrak{v}^{\sf E},D(-\nabla_x)^{\top}A^{\top}D(\nabla_x)\mathbf{U}_{(-)})_{\Gamma_\sigma}\\
&=&(\varepsilon^0)^{\top}M\mathfrak{e}_{(-)}^{\sf E}\ ,
\end{array}
\end{equation}
where $\mathfrak{e}^{\sf M}=D(\nabla_x)\mathfrak{v}^{\sf M}(0)$ and $\mathfrak{e}^{\sf E}=D(\nabla_x)\mathfrak{v}^{\sf E}(0)$.
Here, we apply formula \eqref{e2.64}   for the derivatives $\mathbf{G}^0$ of the Green matrix $G$ in \eqref{e2.64}.
We emphasize that
\begin{equation}\label{e3.222}
(f,\mathbf{G}_{(-)}^0)_\Omega+(g,\mathbf{G}_{(-)}^0)_{\Gamma_\sigma}=-(D(\nabla_x)u)_{(-)}(0)=-\varepsilon_{(-)}^0
\end{equation}
because entries of $\mathbf{G}^0$ are given by the derivatives of
columns of the fundamental matrix $G(x,y)$ with respect to the
second argument, and $G_{(-)}$ satisfies the problem
\begin{equation*}
\begin{array}{l}
D(-\nabla_x)^{\top}A(x)^{\top}D(\nabla_x)G_{(-)}(x,y)=\delta(x-y)\mathbb{I}_{(-)},\ x\in\Omega, \\
D(n(x))^{\top}A(x)^{\top}D(\nabla_x)G_{(-)}(x,y)=0,\ x\in\Gamma_\sigma,\ G_{(-)}(x,y)=0,\ x\in\Gamma_u,
\end{array}
\end{equation*}
where $\mathbb{I}_{(-)}=diag\{1,1,1,-1\}$ (cf. problem
\eqref{e2.64}).

By \eqref{e3.20}   and \eqref{e3.22}, the following representation is valid:
\begin{equation}\label{e3.23}
\mathcal{U}(u^h;\Omega(h))-\mathcal{U}(u;\Omega)=\frac{h^3}{2}(\varepsilon^0)^{\top}M(\mathfrak{e}_{(-)}^{\sf M}-\mathfrak{e}_{(-)}^{\sf E})+O(h^{\alpha+5/2}).
\end{equation}
At the first sight, \eqref{e3.23}   looks like \eqref{e0.1}, however this impression is wrong.
\begin{remark}\label{r2.99}
The  decomposition $u=\mathfrak{v}^{\sf M}+\mathfrak{v}^{\sf E}$
is only a  mathematical device in our analysis, since  in a smart
material   it is impossible to distinguish between the strain
columns $\mathfrak{e}^{\sf M}=D(\nabla_x)\mathfrak{v}^{\sf M}(0)$
and $\mathfrak{e}^{\sf E}=D(\nabla_x)\mathfrak{v}^{\sf E}(0)$
generated at the point $\mathcal{O}$ by the external mechanical
loading $\mathfrak{g}^{\sf M}$ and the electrical surface charge
$\mathfrak{g}^{\sf E}$ in \eqref{e3.21}. Surely, one can measure
only the sum $\varepsilon^0=D(\nabla_x)u(0)$ resulting from
complete  external action and standing as the first term on the
right hand-side of \eqref{e3.23}.
This unususal property of the energy functional should always be taken into account and ignoring the above observation on the decomposition
$u=\mathfrak{v}^{\sf M}+\mathfrak{v}^{\sf E}$
may provoke for misleading physical conclusions. Similar situation
occurs for example for a crack
in a piezoelectric medium. Applying the  Griffith energy fracture
criterion, in \cite{English} the energy release rate at the crack
tip is expressed in terms of stress intensity factors, i.e., local
characteristics of the elastic/electric state at the tip. In
\cite{na349} a mistake in a calculation (formulas (33.23) and
(34.48) in \cite[pages 296 and 312]{English}; cf. comments in
\cite{na349}) was discovered and a non-local formula for the
energy release rate of type \eqref{e3.23}   was derived
rigourously and justified. The non-local character of the energy
release rate means that the energy functional
$\mathcal{U}(u;\Omega)$ cannot be employed for a fracture
criterion and in the Griffith criterion must involve the electric
enthalpy (cf. \cite{ikeda, English} for an interpretation from the
view point of solid fracture mechanics).
\end{remark}

The difference
\begin{equation*}
\mathfrak{e}^{\sf M}-\mathfrak{e}^{\sf E}=(\mathfrak{g}^{\sf E},\mathbf{G}^0)_{\Gamma_\sigma}-(\mathfrak{g}^{\sf M},\mathbf{G}^0)_{\Gamma_\sigma}
\end{equation*}
ought to be regarded as a global characteristics of the mechanical electric state of the body $\Omega$ and, therefore, formula \eqref{e3.23}   has a different physical meaning compared to \eqref{e0.1}   and \eqref{e3.29}   below.

Let us now compute the increment $\mathcal{E}(u^h;\Omega(h))-\mathcal{E}(u;\Omega)$ of the mechanical enthalpy determined in \eqref{e1.16}   and \eqref{e1.17}.
Returning back to the general case $f\neq0$, we obtain
\begin{equation}\label{e3.24}
\begin{array}{l}
\mathcal{E}(u^h;\Omega(h))=\frac{1}{2}(AD(\nabla_x)u^h,D(\nabla_x)u_{(-)}^h)_{\Omega(h)}-(f,u_{(-)}^h)_{\Omega}-(g,u_{(-)}^h)_\Omega\\
=\frac{1}{2}(D(-\nabla_x)^{\top}AD(\nabla_x)u^h,u_{(-)}^h)_{\Omega(h)}+\frac{1}{2}(D(n)^{\top}AD(\nabla_x)u^h,u_{(-)}^h)_{\Gamma_\sigma}\\
-(f,u_{(-)}^h)_{\Omega}-(g,u_{(-)}^h)_\Omega=-\frac{1}{2}(f,u_{(-)}^h)_{\Omega}-\frac{1}{2}(g,u_{(-)}^h)_{\Gamma_\sigma}.
\end{array}
\end{equation}
As above, we have
\begin{equation}\label{e3.25}
(g,u_{(-)}^h)_{\Gamma_\sigma}=(g,u_{(-)})_{\Gamma_\sigma}+h^3(g,\mathbf{U}_{(-)})_{\Gamma_\sigma}+O(h^{\alpha+5/2}).
\end{equation}
Furthermore, in view of representation \eqref{e2.62}   we derive
\begin{equation}\label{e3.26}
(f,u_{(-)}^h)_{\Omega}=(f,u_{(-)})_{\Omega}+h^3(f,\mathbf{U}_{(-)})_{\Omega}+O(h^{\alpha+5/2}).
\end{equation}
according to inequality \eqref{e3.4}   and the following relations
\begin{eqnarray}
h^3|(f,\mathbf{U_{-}})_{\omega_h}|\leq ch^3\int\limits_{\omega_h}|x|^{-2}dx\leq ch^4\leq ch^{\alpha+5/2},\nonumber \\
h|(f,\chi\widetilde{w}^1)_{\Omega(h)}|\leq ch\int\limits_0^{diam\Omega}(1+\frac{r}{h})^{-3}r^2dr\leq ch^4|\ln h|\leq ch^{\alpha+5/2},\label{e3.27} \\
h^2|(f,\chi\widetilde{w}^2)_{\Omega(h)}|\leq ch^2\int\limits_0^{diam\Omega}(1+\frac{r}{h})^{-2+\delta}r^2dr\leq ch^{4-\delta}\leq ch^{\alpha+5/2}.
\end{eqnarray}
In the estimation \eqref{e3.27}   we have used the formulae \eqref{e2.288}   and \eqref{e2.50}   for $\widetilde{w}^1$ and $\widetilde{w}^2$ together with the demanded inclusions $\alpha\in(1/2,1)$ and $\delta\in(0,1/2)$.

Now formulae \eqref{e3.25}, \eqref{e3.26}   and \eqref{e2.67}, \eqref{e3.222}   convert \eqref{e3.24}   into the form
\begin{equation}\label{e3.28}
\begin{array}{l}
\mathcal{E}(u^h;\Omega(h))-\mathcal{E}(u;\Omega)=-\frac{1}{2}(f,u_{(-)})_{\omega_h}+\\
+\frac{1}{2}h^3f(0)^{\top}mes_3\omega((f,G_{(-)}^{0})_{\Omega}+(g,G_{(-)}^{0})_{\Gamma_\sigma})-\\
-\frac{1}{2}h^3(\varepsilon^0)^{\top}M((f,\mathbf{G}_{(-)}^{0})_{\Omega}+(g,\mathbf{G}_{(-)}^{0})_{\Gamma_\sigma})+O(h^{\alpha+5/2})=\\
-\frac{1}{2}((f,u_{(-)})_{\omega_h}-u_{-}(0)^{\top}f(0)mes_3\omega_h)+\frac{1}{2}h^3(\varepsilon^0)^{\top}M\varepsilon_{(-)}^0+O(h^{\alpha+5/2})\\
=\frac{1}{2}h^3(\varepsilon^0)^{\top}M_{(=)}\varepsilon^0+O(h^{\alpha+5/2}).
\end{array}
\end{equation}
  Here, we have taken into account that, first, $M\varepsilon_{(-)}^0=M_{(=)}\varepsilon^0$ according to the definition of $M_{(=)}$ in \eqref{e2.36}   and, second, $u_{-}(0)^{\top}f(0)mes_3\omega_h=(f,u_{(-)})_{\omega_h}+O(h^{3+\alpha})$ due to the smoothness properties \eqref{e2.9}   and \eqref{e2.10}   of $f$ and $u$.

Let us formulate the result obtained in \eqref{e3.28}.
\begin{theorem}\label{t3.5}
The electrical enthalpy \eqref{e1.16}   admits the asymptotic expansion
\begin{equation}\label{e3.29}
\mathcal{E}(u^h;\Omega(h))=\mathcal{E}(u;\Omega)+\frac{1}{2}h^3(\varepsilon^0)^{\top}M_{(=)}\varepsilon^0+O(h^{\alpha+5/2}),
\end{equation}
where $u^h$ and $u$  imply solutions of the piezoelectricity
problems \eqref{e2.2}-\eqref{e2.5}   and
\eqref{e1.7}-\eqref{e1.9}, respectively,
$\varepsilon^0=D(\nabla_x)u(0)$ is the strain column \eqref{e2.11}
and $M_{(=)}=M_{(=)}(A^0,\omega)$ is the modified polarization
matrix which is a symmetric matrix of size $9\times9$ (see
formulae \eqref{e2.35}, \eqref{e2.36}   and Theorem \ref{t2.8}).
\end{theorem}

Note that in contrast to the energy functional \eqref{e1.12}   the
electrical enthalpy has the topological derivative
\begin{equation}\label{e3.30}
\frac{1}{2}(D(\nabla_x)u(0))^{\top}M_{(=)}(A^0,\omega_h)D(\nabla_x)u(0)
\end{equation}
expressed in terms of local characteristics of the
elastic/electric state in the entire body $\Omega$ and of the
shape of the small void $\omega_h$. Owing to
representation \eqref{e2.355}, we emphasize  that the polarization matrix
\eqref{e2.35}   enjoys the homogeneity property
$M(A^0;\omega_h)=h^3M(A^0;\omega)$ which has been used in the
passage from \eqref{e3.29}   to \eqref{e3.30}.

Notice, that exponent $3$ in the factor $h^3$ is conform with formula \eqref{e2.355} for polarization matrix which contains the volume mes$_3\omega$ of the void.

\subsection{Shape functionals and the adjoint state}
\label{3sec:3}

Possible applications of asymptotic analysis performed in the paper include inverse problems, optimum design and shape optimization. We refer the reader to \cite{allaire, fulman1, fulman2} for numerical results of shape and topology optimization  by an application of the levelset method. In the levelset method the topological derivative of a specific shape functional is employed to detect the regions of the hold-all domain to include voids in order to improve the value of the functional to be optimized. The numerical method turn out to be very efficient in two spatial dimensions compared to the pure levelset strategy. Another application with, it seems,  very high potential for numerical solution are all types of inverse problems to detect imperfections within a geometrical domain on the basis of boundary observations. However, in inverse problems it is required that the data imply the unique identification of the imperfection. This property is unknown, in general, for the strategy which is based on the asymptotic analysis in singularly perturbed domains. In particular, it is an open problem how to identify an imperfection from the observation of a finite number of eigenmodes (eigenvalues), which seems to be a natural and efficient way to solve the problem. The difficulty of such an approach is hidden in the fine properties of topological derivatives which are still to be investigated, for example that some sufficiently large set of observations in the mathematical model based on the asymptotic analysis leads to the uniqueness of the position of imperfection. It means that the derivation of topological derivatives is far from being sufficient for the practical applications of the promising tool of shape and topology optimization and identification. We can consider below some specific examples of shape functionals, the other possibilities, including the spectral problem require some additional work to derive the asymptotic formulae.

We return to the analysis.
Recalling the Sobolev embedding theorem $H^1(\Omega)\subset L^6(\Omega)$ in $\mathbb{R}^3$, we assume that the density $J$ in the shape functional
\begin{equation}\label{e3.31}
\mathcal{J}(u;\Omega)=\int\limits_\Omega J(u(x);x)dx
\end{equation}
satisfies the following restrictions:
\begin{equation}\label{e3.32}
|J(a;x)|\leq c(1+|a|^t),
\end{equation}
\begin{equation}\label{e3.33}
|J(b;x)-J(a;x)-J'(a;x)^{\top}(b-a)|\leq c|a-b|^2(1+|a|^{t-2}+|b|^{t-2}),
\end{equation}
\begin{equation}\label{e3.34}
|J(b;x)-J(b;0)|\leq c|x|^{\gamma}(1+|b|^t)
\end{equation}
where $x\in\Omega$, $a$ and $b$ are arbitrary columns in $\mathbb{R}^4$, and the vector function $J'$ is subject to the conditions
\begin{equation}\label{e3.35}
|J'(a;x)|\leq c(1+|a|^{t-1}),
\end{equation}
\begin{equation}\label{e3.355}
|J'(a;x)-J'(b;y)|\leq c(|a-b|^\gamma(|a|^{t-\gamma}+|b|^{t-\gamma})+|x-y|^{\gamma}(|a|^t+|b|^t)),
\end{equation}
while
\begin{equation}\label{e3.36}
t\in[2,6),\ \gamma\in(0,1).
\end{equation}
In other words, along with the restrictions on the  growth of $J$ and $J'$, the integrand $J$ is differentiable with respect to the first variable and H\"older continuous with respect to the second variable.
Moreover, $J'$ is H\"older continuous in both arguments.
Inequality \eqref{e3.32}   ensure that functional \eqref{e3.31}   is defined for $u\in H^1(\Omega)^4\subset L^6(\Omega)^4\subset L^t(\Omega)^4$.

\begin{remark}
\label{RemR7}
Simple examples
$$
\displaystyle\int\limits_{\Omega(h)} R(x)|u^{h}(x) - u(x)|^2 dx
\qquad {\rm and} \qquad
\displaystyle\int\limits_{\Omega(h)}R(x)|u^{h}(x)|^2 dx
$$
are related to the least square method and satisfy the above requirements with $t=2$ and $\gamma=1$ for
$R\in C^{1,\alpha}(\Omega)$. In addition, for $g=0$ in the boundary conditions \eqref{e1.8} and \eqref{e2.3}, the work of external forces
$$
 \displaystyle\int\limits_{\Omega(h)} f^{\sf M}(x)^\top u^{h\sf M}(x) dx
+  \displaystyle\int\limits_{\Omega(h)}f^{\sf E}(x)^\top u^{h\sf E}(x) dx
$$
and the electric enthalpy (cf. \eqref{e3.24})
$$
\mathcal{E}(u^h;\Omega(h))=
\displaystyle-\frac{1}{2}\displaystyle\int\limits_{\Omega(h)} f^{\sf M}(x)^\top u^{h\sf M}(x) dx
+\displaystyle\frac{1}{2}\displaystyle\int\limits_{\Omega(h)}f^{\sf E}(x)^\top u^{h\sf E}(x) dx
$$
readily display another examples.
\end{remark}

We consider the difference
\begin{equation}\label{e3.37}
\mathcal{J}(u^h;\Omega(h))-\mathcal{J}(u;\Omega)=\int\limits_{\Omega(h)}(J(u^{h}(x);x)-J(u(x);x))dx+\int\limits_{\omega_h}J(u(x);x)dx.
\end{equation}
and, owing to \eqref{e3.33}   and \eqref{e2.62}, obtain the formula
\begin{equation}\label{e3.39}
\begin{array}{l}
|J(u^h(x);x)-J(u(x);x)-J'(u(x);x)^{\top}(h^3\mathbf{U}(x)+\chi(x)\sum_{j=1}^{2}h^j\widetilde{w}^j(\frac{x}{h})+\widetilde{u}^h(x))|\leq \\
\qquad \leq
c(h^6|\mathbf{U}(x)|^2+\chi(x)^2\sum_{j=1}^{2}h^{2j}|\widetilde{w}^j(\frac{x}{h})|^2+|\widetilde{u}^h(x)|^2)(1+|u^h(x)|^{t-2}+|u(x)|^{t-2}).
\end{array}
\end{equation}
Recalling the estimates \eqref{e3.4}, \eqref{e3.35}   and applying the H\"older inequality with the index couples $(p,q)=(5/6,6)$ and $(p,q)=(3,2/3)$, we obtain
\begin{equation*}
\begin{array}{l}
\displaystyle\int\limits_{\Omega(h)}J'(u(x);x)^{\top}\widetilde{u}^h(x)dx\leq c\displaystyle\int\limits_{\Omega(h)}(1+|u(x)|^5)|\widetilde{u}^h(x)|dx\leq \\
\qquad\qquad\leq
c(1+\|u;L^6(\Omega)\|^5)\|\widetilde{u}^h;L^6(\Omega(h))\|\leq
c\|\widetilde{u}^h;H^1(\Omega)\|\leq
ch^{\alpha+5/2},\\
\displaystyle\int\limits_{\Omega(h)}|\widetilde{u}^h|^2(1+|u^h|^{t-2}+|u|^{t-2})dx\leq c\displaystyle\int\limits_{\Omega(h)}|\widetilde{u}^h|^2(1+|u^h|^{4}+|u|^{4})dx\leq\\
\leq c\|\widetilde{u}^h;L^6(\Omega(h))\|^2(1+\|u^h;L^6(\Omega(h))\|^4+\|u;L^6(\Omega)\|^4)\leq ch^{2\alpha+5}.
\end{array}
\end{equation*}
Similarly,
\begin{equation*}
h^6\int\limits_{\Omega\setminus\mathbb{B}_R'}|\mathbf{U}(x)|^2(1+|u^h(x)|^{t-2}+|u(x)|^{t-2})dx\leq ch^6.
\end{equation*}
However, because of the singularity $|\mathbf{U}(x)|=O(|x|^{-2})$,  we use in the ball $\mathbb{B}_{R'}$ the H\"older inequality with the couple
\begin{equation}\label{e3.38}
(p,q)=\left(\frac{6}{8-t},\frac{6}{t-2}\right)
\end{equation}
to derive that
\begin{equation*}
\begin{array}{l}
h^6\displaystyle\int\limits_{\mathbb{B}_{R'}\setminus\omega_h}|\mathbf{U}|^2(1+|u^h|^{t-2}+|u|^{t-2})dx\leq ch^6\left(\displaystyle\int\limits_{ch}^{R'}
\displaystyle {r^{-\frac{24}{8-t}}}r^2dr\right)^{\frac{8-t}{6}}\times\\
\qquad\qquad\times(1+\|u^h;H^1(\Omega(h))\|^{t-2}+\|u;H^1(\Omega)\|^{t-2})\leq
ch^{6-t/2}.
\end{array}
\end{equation*}

We deal with the boundary layers in the same way as in
\eqref{e3.15}   and \eqref{e3.155}. Outside the ball
$\mathbb{B}_{Rh}$ we apply the inequalities \eqref{e2.288}   and
\eqref{e2.50}   even much rougher  ones, to conclude by the
H\"older inequality with the index couple \eqref{e3.38}   that
\begin{equation}\label{e3.40}
\begin{array}{l}
h^{2j}\displaystyle\int\limits_{\Omega\setminus\mathbb{B}_{Rh}}\left|\chi(x)\widetilde{w}^j\left(\frac{x}{h}\right)\right|^2(1+|\widetilde{u}^h(x)|^2+|u(x)|^2)dx\leq\\
\qquad\qquad\leq
ch^6\left(\displaystyle\int\limits_{ch}^{R'}r^{-\frac{12(3-j)}{8-t}}r^2dr\right)^{\frac{8-t}{6}}\leq
ch^{6-t/2},\ j=1,2.
\end{array}
\end{equation}
Inside the ball $\mathbb{B}_{Rh}$ the H\"older inequality gives
\begin{equation*}
\begin{array}{l}
h^{2j}\displaystyle\int\limits_{\mathbb{B}_{Rh}\setminus\omega_h}\left|\widetilde{w}^j\left(\frac{x}{h}\right)\right|^2(1+|\widetilde{u}^h(x)|^2+|u(x)|^2)dx\\
\qquad\qquad\leq ch^{2j}\left(\displaystyle\int\limits_{\mathbb{B}_{Rh}\setminus\omega_h}|\widetilde{w}^j\left(\frac{x}{h}\right)|^{\frac{12}{8-t}}dx\right)^{\frac{8-t}{6}}=ch^{2j+3\frac{8-t}{6}}\left(\displaystyle\int\limits_{\mathbb{B}_{Rh}\setminus\omega}\widetilde{w}^j\left(\frac{x}{h}\right)|^{\frac{12}{8-t}}d\xi\right)^{\frac{8-t}{6}}\\
\qquad\qquad\leq ch^{2j+4-t/2}\leq ch^{6-t/2},\ j=1,2.
\end{array}
\end{equation*}
Note that $\frac{12}{8-t}<6$ due to \eqref{e3.36}   and, therefore,
\begin{equation*}
\|\widetilde{w}^j;L^{\frac{12}{8-t}}(\mathbb{B}_R\setminus\omega)\|\leq c\|\widetilde{w}^j;H^1(\mathbb{B}_R\setminus\omega)\|\leq c\|\widetilde{w}^j;V_0^1(\Xi)\|.
\end{equation*}

Although,  the faster rates of decay of the remainders
$\widetilde{w}^1$ and $\widetilde{w}^2$ (cf. \eqref{e2.55}) are
not used in the estimation \eqref{e3.40}, the rate of decay
becomes an important ingredient of the inequalities
\begin{equation*}
h^j|\int\limits_{\Omega(h)}J'(u(x);x)^{\top}\chi(x)\widetilde{w}^j(\frac{x}{h})dx|\leq ch^{7/2},\ j=1,2,
\end{equation*}
its  derivation is much simpler,  though.
A simplification originates from the relation $|J'(u(x);x)|\leq
const$ for $x\in supp\chi\subset\mathbb{B}_{R'}$ so that one may
repeat the calculation \eqref{e3.27}.

Finally, we write
\begin{equation*}
h^3|\int\limits_{\omega_h}J'(u(x);x)^{\top}\mathbf{U}(x)dx|\leq ch^3\int\limits_{0}^{Rh}r^{-2}r^2dr\leq ch^4
\end{equation*}
and, in view of \eqref{e2.10}   and \eqref{e3.355},
\begin{equation*}
\left|\int\limits_{\omega_h}J(u(x);x)dx-h^3J(u(0);0)mes_3\omega\right|\leq ch^{3+min\{\alpha,\gamma\}}
\end{equation*}

Everything is prepared to derive a formula of type \eqref{e0.1}   for the  shape functional \eqref{e3.31}.
\begin{theorem}\label{t3.81}
Let the assumption formulated above hold true.
Then the asymptotic formula
\begin{equation}\label{e3.41}
\begin{array}{rcl}
\mathcal{J}(u^h;\Omega(h))&=&\mathcal{J}(u;\Omega)
+ h^3((J(u(0);0)-P(0)^{\top}f(0))mes_3\omega
\\
&-&(D(\nabla_x)P(0))^{\top}M\varepsilon^0)+O(h^{3+min\{\gamma,\alpha-1/2,3-t/2\}})
\end{array}
\end{equation}
is valid where $P\in\mathring{H}^1(\Omega;\Gamma_u)^4\cap C^{2,min\{\alpha,\gamma\}}(\mathbb{B}_{R'})^4$ is a solution of the formally adjoint piezoelectricity problem
\begin{equation}\label{e3.42}
\begin{array}{l}
D(-\nabla_x)^{\top}A(x)^{\top}D(\nabla_x)P(x)=J'(u(x);x),\ x\in\Omega,\\
D(n(x))^{\top}A(x)^{\top}D(\nabla_x)P(x)=0,\ x\in\Gamma_\sigma, P(x)=0,\ x\in\Gamma_\sigma.
\end{array}
\end{equation}
\end{theorem}

{\bf Proof.}
The calculations performed above provide the relation
\begin{equation*}
h^{-3}(\mathcal{J}(u^h;\Omega(h))-\mathcal{J}(u;\Omega))=J(u(0);0)mes_3\omega+(J'(u),\mathbf{U})_\Omega+O(h^{min\{\gamma,\alpha-1/2,3-t/2\}}).
\end{equation*}
We recall the representation \eqref{e2.67}   where $G^0$ is the
Green matrix, i.e., a solution to the problem \eqref{e2.64}. The
Green matrix and its derivatives help to calculate the solution $P$ of
the formally adjoint problem \eqref{e3.42}   and the derivatives
$\mathbf{G}^0$ (see \eqref{e2.65}) deliver the column
$D(\nabla_x)P(x)$ at the point $x=0$. In other  words, we write
\begin{equation}\label{e3.44}
\begin{array}{rcl}
(J'(u),\mathbf{U})_\Omega&=&(D(-\nabla_x)^{\top}A^{\top}D(\nabla_x)P,\mathbf{G}^0)_\Omega M\varepsilon^0
\\
&-&mes_3 \omega (D(-\nabla_x)^{\top}A^{\top}D(\nabla_x)P, {G}^0)_\Omega f(0) \\
&=&(P,D(\nabla_x)^{\top}\delta M\varepsilon^0)_\Omega-mes_3\omega(P,\delta f(0))_\Omega\\
&=&-(D(\nabla_x)P(0))^{\top}M\varepsilon^0-P(0)^{\top}f(0)mes_3\omega.
\end{array}
\end{equation}
We again used the Dirac mass $\delta$ in the framework of the
theory of distributions to compute the expression \eqref{e3.44}.

Finally, in order to justify our calculations we make the following comments.
By assumptions \eqref{e3.32}, \eqref{e3.36}   and \eqref{e3.355}, \eqref{e2.8}, the functional
\begin{equation*}
\mathring{H}^1(\Omega;\Gamma_u)^4\ni v\rightarrow(J'(u),v)_\Omega
\end{equation*}
is continuous and $J'(u)\in
C^{0,min\{\alpha,\gamma\}}(\mathbb{B}_{R'})$ with any $R'<R$.
Thus, the same arguments as in Sections \ref{3sec:1} and
\ref{1sec:2} guarantee the existence of a solution $P$ to the
problem \eqref{e3.42}   which is twice differentiable in the
vicinity of the point $x=0$. These observations make all
calculations justified.$\blacksquare$

The topological derivative of the functional $\mathcal{J}$, i.e.,
\begin{equation*}
\mathcal{T}(u,\omega)=(J(u(0);0)-P(0)^{\top}f(0))mes_3\omega-(D(\nabla_x)P(0))^{\top}M\varepsilon^0,
\end{equation*}
is non-local since it involves the adjoint state $P$ in
\eqref{e3.42}   which depends on the solution $u$ of the
piezoelectricity problem in the entire domain $\Omega$.

\subsection{Example}
\label{sec:examples} It turns out, that the so-called weak
interaction is quite common feature of piezoelectric materials.
Therefore, we are going to present an example for such materials.
Assume that there is a weak interaction between the mechanical and
electric fields. This means that in the decomposition

\begin{gather}
\label{e5.1}
A=A_{(0)} + A_{(1)}\ ,
\\
\notag
A_{(0)}=\left(
\begin{array}{cc}
A^{\sf{MM}} &\mathbb{O}_{6 \times 3} \\
\mathbb{O}_{3 \times 6} & A^{\sf{EE}}
\end{array}
\right)
\ ,\quad
A_{(1)}=\left(
\begin{array}{cc}
\mathbb{O}_{6 \times 6} & -A^{\sf{ME}}\\
A^{\sf{EM}} & \mathbb{O}_{3 \times 3}
\end{array}
\right)
\end{gather}
the entries of matrix $ A_{(1)}$ are much smaller compared to non
trivial entries of the matrix $A_{(0)}$. It implies that in the
first order approximation the piezoelectricity problem is
decoupled into two problems, the pure elasticity problem with the
stiffness matrix $A^{\sf{MM}}$, and the pure electricity problem
with the permeability matrix $A^{\sf{EE}}$.

We are going to evaluate the main correction terms in the
asymptotic expansions of characteristics for the piezoelectric
bodies $\Omega$, $\Xi$ and $\Omega(h)$ (see Sections \ref{2sec:1},
\ref{3sec:2} and \ref{1sec:2}).  We point out that to evaluate the main asymptotic terms in all formulae given below, it is sufficient to solve only the pure elasticity and the pure electricity problems.

\begin{remark}
\label{dots}
Since we always deal with the first order asymptotic corrections, the introduction of any small amplitude parameter neither makes formulae more transparent, nor contribute to the exactness of presentation. We emphasize that, in contrast to the preceding sections, the perturbations here are of regular type, which means that the justification of obtained formulae relies upon the standard argument of convergence of Neumann series. In order to simplify the notation, in the sequel the second order terms are always denoted by dots, starting from \eqref{e5.2}.
\end{remark}

We proceed with the solution
\begin{gather}
\label{e5.2}
u(x)=u_{(0)} + u_{(1)} (x)+ \dots
\end{gather}
of the problem  \eqref{e1.7}-\eqref{e1.9}.  In view of
\eqref{e5.1}, the displacement vector $u_{(0)}^{\sf{M}}$ and the
electric vector $u_{(0)}^{\sf{E}}$ verify the problems

\begin{gather}
\label{e5.3}
D^{\sf{M}}(-\nabla_x)^{\top}A^{\sf{MM}}(x)D^{\sf{M}}(\nabla_x)u_{(0)}^{\sf{M}}(x)=f^{\sf{M}}(x),\ x\in\Omega,
\\
\notag
D^{\sf{M}}(n(x))^{\top}A^{\sf{MM}}(x)D^{\sf{M}}(\nabla_x)u_{(0)}^{\sf{M}}(x)=g^{\sf{M}}(x),\ x\in\Gamma_{\sigma},\quad u_{(0)}^{\sf{M}}(x)=0,\
x\in\Gamma_{u},\\
\label{e5.4}
-\nabla_x^{\top} A^{\sf{EE}} \nabla_x u_{(0)}^{\sf{E}}(x) = f^{\sf{E}}(x),\ x\in\Omega,
\\
\notag
n^{\top} A^{\sf{EE}}  u_{(0)}^{\sf{E}}(x) = g^{\sf{E}}(x),\ x\in\Gamma_{\sigma},\quad  u_{(0)}^{\sf{E}}(x)=0,\ x\in\Gamma_{u}\ ,
\end{gather}
and can be determined separately. Inserting  \eqref{e5.2} and  \eqref{e5.1} into  \eqref{e1.7}- \eqref{e1.9}, we arrive at the problem

\begin{gather}
\notag
D(-\nabla_x)^{\top}A_{(0)}(x)D(\nabla_x)u_{(1)}(x)=D(-\nabla_x)^{\top}A_{(1)}(x)D(\nabla_x)u_{(0)}(x),\ x\in\Omega, \\
D(n(x))^{\top}A_{(0)}(x)D(\nabla_x)u_{(1)}(x)=D(n(x))^{\top}A_{(1)}(x)D(\nabla_x)u_{(0)}(x),\ x\in\Gamma_{\sigma},\label{e5.5} \\
\notag
u_{(1)}(x)=0,\ x\in\Gamma_{u}\ .
\end{gather}
This problem is decoupled as well, however,  its solution manifests the interaction between electric and mechanical fields, since the displacement vector $u^{\sf{M}}_{(1)}$ depends only on the main part $u^{\sf{E}}_{(0)}$ of the electric potential and, in the same manner, $u^{\sf{E}}_{(1)}$ depends on $u^{\sf{M}}_{(0)}$.

In order to complete the asymptotic formulae, in the same way as in the previous sections, we also need the expansion for the polarization matrix

\begin{gather}
\label{e5.6}
M=M_{(0)} + M_{(1)} + \dots \ ,
\\
\notag
M_{(0)}=\left(
\begin{array}{cc}
M_{(0)}^{\sf{M}} &\mathbb{O}_{6 \times 3} \\
\mathbb{O}_{3 \times 6} & M_{(0)}^{\sf{E}}
\end{array}
\right)
\ ,\quad
M_{(1)}=\left(
\begin{array}{cc}
\mathbb{O}_{6 \times 6} & M_{(1)}^{\sf{ME}}\\
M_{(1)}^{\sf{EM}} & \mathbb{O}_{3 \times 3}
\end{array}
\right)\ .
\end{gather}

We emphasize that the matrices $M_{(0)}$ and $M_{(1)}$ inherit the block diagonal structure of $A_{(0)}$ and the block-anti-diagonal of $A_{(1)}$, respectively. The same structures are kept by all matrix objects, in particular, the fundamental matrix takes the form

\begin{gather}
\label{e5.7}
\Phi =\Phi _{(0)} + \Phi _{(1)} + \dots \ ,
\\
\notag
\Phi _{(0)}=\left(
\begin{array}{cc}
\Phi _{(0)}^{\sf{M}} &\mathbb{O}_{3 \times 1} \\
\mathbb{O}_{1 \times 3} & \Phi _{(0)}^{\sf{E}}
\end{array}
\right)
\ ,\quad
\Phi _{(1)}=\left(
\begin{array}{cc}
\mathbb{O}_{3 \times 3} & \Phi _{(1)}^{\sf{ME}}\\
\Phi _{(1)}^{\sf{EM}} & 0
\end{array}
\right)\ .
\end{gather}
Here, $\Phi _{(0)}^{\sf{M}} $ is the fundamental matrix for the
elasticity matrix operator
$D^{\sf{M}}(-\nabla_\xi)^{\top}A^{0\sf{M}}D^{\sf{M}}(\nabla_\xi)$
and $\Phi _{(0)}^{\sf{E}}$ is the fundamental matrix for the scalar
operator $-\nabla_\xi^{\top} A^{0\sf{E}} \nabla_\xi$. Furthermore,
$M_{(0)}^{\sf{E}} $ and $M_{(0)}^{\sf{M}} $ are the virtual mass matrix and
the elasticity polarization matrix for the cavity
$\omega\subset\mathbb{R}^3$, which are negative definite (see
\cite{PoSe} and \cite{na280, NazSokSPN}).

It is convenient to proceed with the matrix solution \eqref{e2.33}
which, according to  \eqref{e5.1} and \eqref{e2.34}, enjoys the
expansion

\begin{gather}
\label{e5.8}
W=W_{(0)} + W_{(1)} + \dots\ ,
\\
\notag
W_{(0)}=\left(
\begin{array}{cc}
W_{(0)}^{\sf{M}} &\mathbb{O}_{3 \times 3} \\
\mathbb{O}_{1 \times 6} & W_{(0)}^{\sf{E}}
\end{array}
\right)
\ ,\quad
W_{(1)}=\left(
\begin{array}{cc}
\mathbb{O}_{3 \times 6} & W_{(1)}^{\sf{ME}}\\
W_{(1)}^{\sf{EM}} & \mathbb{O}_{1 \times 3}
\end{array}
\right)
\end{gather}
with

\begin{gather}
\label{e5.9}
W(\xi)=(M D(\nabla_\xi) \Phi(\xi)^\top)^\top + O(|\xi|^{-2})=\\
\notag
(M _{(0)} D(\nabla_\xi) \Phi _{(0)}(\xi)^\top)^\top + (M _{(0)} D(\nabla_\xi )\Phi _{(1)}(\xi)^\top + M _{(1)} D(\nabla_\xi \Phi _{(0)}(\xi)^\top)^\top +\dots +O(|\xi|^{-2})\ .
\end{gather}

The correction term $\Phi _{(1)}$ in \eqref{e5.7} is a power-law solution of form \eqref{e2.45} for the system of differential equations

\begin{gather}
\label{e5.10}
D(-\nabla_\xi)^\top A^0 _{(0)}D(\nabla_\xi)\Phi _{(1)}(\xi)
=
D(\nabla_\xi)^\top A^0 _{(1)}D(\nabla_\xi)\Phi _{(0)}(\xi)\ ,\quad \xi\in \mathbb{R}^3\setminus \{0\}\ ,
\end{gather}
(cf. \eqref{e2.46}).  By a general result in \cite{Ko} (see also
\cite[Lemmas 3.3.1 and 3.5.11]{NaPl}), the solution $\Phi _{(1)}$
can depend linearly on $\ln{|\xi |}$, however, the same argument
as in the proof of Lemma \ref{l2.9} ensures that $\Phi _{(1)}$ is
positive homogeneous of degree -$1$ according to \eqref{e2.27}.
The solution $\Phi _{(1)}$, which  is defined up to the linear
combination $\Phi _{(0)}C$ of the fundamental matrix columns with
the constant column $C\in \mathbb{R}^4$, can be fixed such that
\begin{gather}
\label{e5.11}
\int\limits_{\mathbb{S_1}}D(\nabla_\xi)^\top A^0 D(\nabla_\xi)\Phi _{(1)}(\xi) ds_\xi =0\in \mathbb{R}^4\ .
\end{gather}

The exterior problem for the correction term in \eqref{e5.8} takes
the form

\begin{gather}
\label{e5.12}
D(-\nabla_\xi)^\top A^0 _{(0)}D(\nabla_\xi) W_{(1)}(\xi) =
D(\nabla_\xi)^\top A^0 _{(1)}D(\nabla_\xi) W_{(0)}(\xi)\ ,\xi\in\Xi\ ,\\
\label{e5.13}
D(n^\omega(\xi))^\top A^0 _{(0)}D(\nabla_\xi) W_{(1)}(\xi) =
-D(n^\omega(\xi))^\top A^0 _{(1)}D(\nabla_\xi) W_{(0)}(\xi)\ ,\xi\in\partial\omega\ .
\end{gather}
Since, owing to  \eqref{e5.9}, we have
\begin{gather*}
W_{(0)}(\xi) =(M _{(0)}D(\nabla_\xi)\Phi _{(0)}(\xi)^\top )^\top
+ O(|\xi|^{-3})\ ,
\end{gather*}
the right-hand side $F_{(1)}(\xi)$ in \eqref{e5.12} admits the decomposition

\begin{gather}
\notag
F_{(1)}(\xi)=D(\nabla_\xi)^\top A^0 _{(1)}D(\nabla_\xi) (D(\nabla_\xi)\Phi _{(0)}(\xi)^\top )^\top M _{(0)}^\top
+\widetilde F_{(1)}(\xi)=
\\
\label{e5.14}
=\sum_{q=1}^3 \frac{\partial}{\partial \xi_q}D(\nabla_\xi)^\top A^0 _{(1)}D(\nabla_\xi) \Phi _{(0)}(\xi)^\top D(e_q)^\top M _{(0)}^\top+\widetilde F_{(1)}(\xi)=
\\
\notag
=\sum_{q=1}^3 D(-\nabla_\xi)^\top A^0 _{(0)}D(\nabla_\xi)\frac{\partial \Phi _{(1)}}{\partial \xi_q}D(e_q)^\top M _{(0)}^\top+\widetilde F_{(1)}(\xi)
\end{gather}
with the remainder $\widetilde F_{(1)}(\xi)=O(|\xi|^{-5})$. In
\eqref{e5.14}, the equation \eqref{e5.10} has been applied.
Comparing \eqref{e5.14} with \eqref{e5.9}, we set

\begin{gather}
\label{e5.15}
W_{(1)}(\xi)=\widetilde W_{(1)}(\xi) + (M _{(0)}D(\nabla_\xi) \Phi _{(0)}(\xi)^\top)^\top \ .
\end{gather}
Recall that $\omega$ contains the origin $\xi=0$, therefore, the last term in  \eqref{e5.15} is smooth in $\overline\Xi$. As a result, a new exterior problem is obtained, with the right-hand side $\widetilde F_{(1)}$ which decays sufficiently fast at infinity,

\begin{gather}
\label{e5.16}
D(-\nabla_\xi)^\top A^0 _{(0)}D(\nabla_\xi) \widehat W_{(1)}(\xi) =\widetilde F_{(1)}(\xi)
\ ,\xi\in\Xi\ ,\\
\notag
D(n^\omega(\xi))^\top A^0 _{(0)}D(\nabla_\xi) \widehat W_{(1)}(\xi) =\widetilde G_{(1)}(\xi)
\ ,\xi\in\partial\omega\ ,
\end{gather}
where
\begin{gather}
\label{e5.1600}
\widetilde G_{(1)}(\xi)=
D(n^\omega(\xi))^\top A^0 _{(1)}D(\nabla_\xi) \widehat W_{(0)}(\xi)-\\
\notag
- D(n^\omega(\xi))^\top A^0 _{(0)}D(\nabla_\xi)  (M _{(0)}D(\nabla_\xi) \Phi _{(1)}(\xi)^\top)^\top \ .
\end{gather}
Now, the decay of $\widetilde G_{(1)}(\xi)$  can be used, indeed,
by Proposition \ref{p2.4} (see \cite{Ko} and \cite[Theorem
3.5.6]{NaPl}) and the calculations  \eqref{e2.30}, \eqref{e2.29},
the solution $\widehat W_{(1)}\in V^1_0(\Xi)^4$ admits the
asymptotic form
\begin{gather}
\label{e5.18}
\widehat W_{(1)}(\xi)= (M _{(1)}D(\nabla_\xi) \Phi _{(0)}(\xi)^\top)^\top + \widetilde W_{(1)}(\xi)\ ,
\end{gather}
where the remainder $ \widetilde W_{(1)}$ is subject to the
estimates \eqref{e2.288} with the majorants $c_k
\rho^{-3-k+\delta}$ ($\delta>0$ is arbitrary) and the notation
used for the derivatives of the fundamental matrix $\Phi _{(0)}$
is matched with formulae  \eqref{e5.9} and \eqref{e5.15}.

In order to evaluate the correction term $M_{(1)}$ in  the
expansion of the polarization matrix the method \cite{MaPl1} is
employed, here we recall that the columns of the matrix

\begin{gather}
\label{e5.19}
\mathcal{W}_{(0)(-)}(\xi)=D_{(-)}(\xi)^\top +  {W}_{(0)(-)}(\xi)
\end{gather}
(cf.  \eqref{e2.355})   are formal solutions to the homogeneous problem \eqref{e5.16}. By the Green formula in $\Xi\cap\mathbb{B}_R$, we obtain

\begin{gather}
\label{e5.20}
\int\limits_{\Xi\cap \mathbb{B}_R}
\mathcal{W}_{(0)(-)}(\xi)^\top \widetilde F_{(1)}(\xi) d\xi +
\int\limits_{\partial\omega} \mathcal{W}_{(0)(-)}(\xi)^\top \widetilde G_{(1)}(\xi) ds_\xi=\\
\notag
\int\limits_{\partial \mathbb{B}_R}
( \widehat W_{(1)}(\xi)^\top D(|\xi|^{-1}\xi)^\top A^0 _{(0)(-)}D(\xi)\mathcal{W} _{(0)(-)}-
\mathcal{W} _{(0)(-)}^\top  D(|\xi|^{-1}\xi)^\top A^0 _{(0)(-)}D(\xi) \widehat W_{(1)}(\xi)) ds_\xi
\\
\notag
+ O(R^{-1})= -M_{(1)(=)} + O(R^{-1})\ .
\end{gather}

We have here repeated the computation \eqref{e2.38} based on the
representations \eqref{e5.18} and  \eqref{e5.19}. The integrand on
the left-hand side of \eqref{e5.20} is of order $|\xi|^{-4}$ and,
hence, the integral over $\Xi$ converges and the formula

\begin{gather}
\label{e5.21}
M_{(1)(=)}
=\left(
\begin{array}{cc}
\mathbb{O}_{6 \times 6} & M_{(1)}^{\sf{ME}}\\
M_{(1)}^{\sf{EM}} & \mathbb{O}_{3 \times 3}
\end{array}
\right)
 =- \int\limits_{\Xi} \mathcal{W}_{(0)(-)}(\xi)^\top \widetilde F_{(1)}(\xi) d\xi +
\int\limits_{\partial\omega} \mathcal{W}_{(0)(-)}(\xi)^\top \widetilde G_{(1)}(\xi) ds_\xi
\end{gather}
together with \eqref{e5.14}-\eqref{e5.16} expresses the matrix
$M_{(1)}$ (cf. the definition \eqref{e2.36}) in terms of the
matrix $A^{\sf{ME}}=(A^{\sf{EM}})^{-1}$ and the special solutions
$W^1,\dots,W^6 $ and $W^7,W^8,W^9 $ of the pure elasticity and the
pure electricity exterior problems in $\Xi$. Theorem \ref{t2.8}
shows that $(M^{\sf{ME}}_{(1)})^\top =-M^{\sf{EM}}_{(1)}$.

The formulae derived above can be used, e.g., to obtain the
topological derivative of the electric enthalpy \eqref{e3.29}:

\begin{gather}
\label{e5.22}
\mathcal{T}_{\mathcal{E}}(u;\omega)=
\\
\notag
=\frac{1}{2} h^3
(( D^{\sf{M}}(\nabla_x) u^{\sf{M}}_{(0)}(0))^\top M_{(0)}^{\sf{M}} D ^{\sf{M}}
(\nabla_x) u^{\sf{M}}_{(0)}(0)-
\nabla_x u^{\sf{E}}_{(0)}(0)^\top M_{(0)}^{\sf{E}}
\nabla_x u^{\sf{E}}_{(0)}(0)) +
\\
\notag
+h^3
(( D^{\sf{M}}(\nabla_x) u^{\sf{M}}_{(0)}(0))^\top M_{(0)}^{\sf{M}} D ^{\sf{M}}
(\nabla_x) u^{\sf{M}}_{(1)}(0)-
\nabla_x u^{\sf{E}}_{(0)}(0)^\top M_{(0)}^{\sf{E}}
\nabla_x u^{\sf{E}}_{(1)}(0)) +
\\
\notag
+ h^3
\nabla_x u^{\sf{E}}_{(0)}(0))^\top M^{\sf{EM}}
D ^{\sf{M}}(\nabla_x) u^{\sf{M}}_{(1)}(0)+\dots,
\end{gather}
where $M_{(0)}^{\sf{M}}$ and $M_{(0)}^{\sf{E}}$ are the elasticity polarization matrix and the virtual mass matrix for the cavity while $M_{(1)}^{\sf{EM}}=-\left( M_{(1)}^{\sf{ME}}\right)^\top$ is expressed in \eqref{e5.21}.

Even the main term (with the factor $\frac{1}{2} h^3$) of the
topological derivative \eqref{e5.22} has no sign, that is, in
contrast to the forms of topological derivatives of the energy
functionals for the pure elasticity and the pure electricity
problems. The correction term (with factor $h^3$) in \eqref{e5.22}
depends on two specific ingredients, namely, the correction term
$M^{\sf{EM}}$ in polarization matrix (see  \eqref{e5.6}  and
\eqref{e5.21}),   and the correction terms $u^{\sf{M}}_{(1)},
u^{\sf{E}}_{(1)}$ for the combined mechanical and electric fields.

\begin{remark}
\label{r01.000} All the attributes in the above formulae  can be
given explicitly for some canonical shapes, including balls, ellipsoids and elliptic cracks in
three spatial dimensions, and some other shapes in two spatial
dimensions (see \cite{PoSe} and  \cite{na100, MovMov, pos2, ammari} and others).
\end{remark}

\begin{remark}
\label{r00.000} The case of  $g^{\sf{E}}=0, f^{\sf{E}}=0$ has a
very clear physical meaning (i.e. one gets an electric sparkle
when pressing the lighter button). Then, in notation of Section
\ref{2sec:3},
\begin{gather*}
\mathfrak{u}^{\sf{M}}=u,\quad \mathfrak{u}^{\sf{E}}=0,\quad \mathfrak{e}^{\sf{M}}=\varepsilon^0, \quad\mathfrak{e}^{\sf{E}}=0,
\end{gather*}
thus, by relation  \eqref{e2.36}, we can conclude that the
topological derivatives in \eqref{e2.36} and \eqref{e2.39} of the
energy and electric enthalpy functionals coincides one with
another. In general, this identity is false, and can be misleading
for the choice of governing Gibbs' functional for piezoelectric
body (cf. Remark  \ref{r2.99}). The relations between the
topological derivatives for elasticity and piezoelectricity are
easy to established, since the topological derivative for
piezoelectricity can be viewed as the difference of that for
elasticity and of the other for electricity.
\end{remark}

\bigskip

\textbf{Acknowledgements.} This paper was prepared during the visit of S.A. Nazarov to the
Institute Elie Cartan of the University Henri Poincar\'e Nancy 1
and to Department of Civil Engineering of Second University of
Naples. The research of S.A.N. is partially supported by the grant
RFFI-09-01-00759 and by project "Asymptotic analysis of composite
materials and thin and non-homogeneous structures" (Regione
Campania, law n.5/2006); J.S. was partially supported by the
Projet CPER Lorraine MISN: Analyse, optimisation et contr\^ole 3
in France, and the grant N51402132/3135 Ministerstwo Nauki i
Szkolnictwa Wyzszego: Optymalizacja z wykorzystaniem pochodnej
topologicznej dla przeplywow w osrodkach scisliwych in Poland.

\end{document}